\documentclass[10pt]{article}
\usepackage{amsmath}
\usepackage{float}
\usepackage{amsfonts}
\usepackage{amssymb}
\usepackage{amsmath}
\usepackage{tikz}
\usepackage{paralist}
\usepackage{epstopdf}
\usepackage[a4paper, total={5.5in, 7.5in}]{geometry}
\usepackage{hyperref}

\hypersetup{
colorlinks=true,
linkcolor=blue,
filecolor=magenta,      
urlcolor=cyan,
}
\hypersetup{
colorlinks=true,
linkcolor=blue,
filecolor=magenta,      
urlcolor=cyan,
}
\newtheorem{theorem}{Theorem}

\newtheorem{definition}[theorem]{Definition}
\newtheorem{example}[theorem]{Example}

\newtheorem{lemma}[theorem]{Lemma}

\newtheorem{proposition}[theorem]{Proposition}
\newtheorem{remark}[theorem]{Remark}

\newenvironment{proof}[1][Proof]{\noindent\textbf{#1.} }{\ \rule{0.5em}{0.5em}}

\begin{document}

\title{On boundary controllability and stabilizability of the 1D wave
equation in non-cylindrical domains}
\author{Mokhtari Yacine}
\maketitle
\tableofcontents

\begin{abstract}
In this paper, we deal with boundary controllability and \ boundary
stabilizability of the 1D wave equation in non-cylindrical domains. By using the characteristics method, we
prove under a natural assumption on the boundary functions that the 1D wave
equation is controllable and stabilizable from one side of the boundary.
Furthermore, the control function and the decay rate of the solution are
given explicitly.
\end{abstract}

\footnotetext{%
date: 02/03/2020} \footnotetext{%
E-mail adress: \href{mailto:yacine.mokhtari@univ-fcomte.fr}{%
yacine.mokhtari@univ-fcomte.fr}} \footnotetext{%
Laboratoire de Math\'{e}matiques UMR 6623, Universit\'{e} de Bourgogne
Franche-Comt\'{e}, 16, route de Gray, 25030 Besan\c{c}on cedex, France} 
\footnotetext{%
University of Sciences and Technology Houari Boumedienne P.O.Box 32, El-Alia
16111, Bab Ezzouar, Algiers, Algeria}

	\section{Introduction and preliminaries}
	
	In this work, we are interested in the boundary controllability and
	stabilizability of the one dimensional wave equation in non-cylindrical
	domains. More precisely, let $\alpha $ and $\beta $ be two real functions
	defined on $%
	\mathbb{R}
	_{+}$ and $Q$ be the set 
	\begin{equation*}
	Q=\left\{ (t,x)\in 
	\mathbb{R}
	^{2},\text{ }x\in (\alpha (t),\beta (t)),\text{ }\alpha (t)<\beta (t),\text{ 
	}t\in (0,\infty )\right\} ,
	\end{equation*}%
	with $\alpha (0)=0$ and $\beta (0)=1$. We consider the following two systems 
	\begin{equation}
	\left\{ 
\begin{array}{ccc}
	y_{tt}(t,x)=y_{xx}(t,x),\text{ \ \ \ \ \ \ \ \ \ \ \ \ \ \ \ \ } & \mathrm{\
		in} & Q,\text{ \ \ \ \ } \\ 
	y(t,\alpha (t))=\frac{1}{2}u(t),\text{ }y(t,\beta (t))=0, & \mathrm{in} & 
	(0,\infty ), \\ 
	y(0,x)=y_{0}(x),\text{ }y_{t}(0,x)=y_{1}(x), & \mathrm{in} & (0,1),\text{ }%
	\end{array}%
	\right.  \label{control}
	\end{equation}%
	and 
	\begin{equation}
	\left\{ 
	\begin{array}{ccc}
	y_{tt}(t,x)=y_{xx}(t,x),\text{ \ \ \ \ \ \ \ \ \ \ \ \ \ \ \ \ \ \ \ \ \ \ \
		\ \ \ \ \ } & \mathrm{in} & Q,\text{ \ \ \ \ } \\ 
	y_{t}(t,\alpha (t))=f(t)y_{x}(t,\alpha (t)),\text{ }y(t,\beta (t))=0, & 
	\mathrm{in} & (0,\infty ), \\ 
	y(0,x)=y_{0}(x),\text{ }y_{t}(0,x)=y_{1}(x),\text{ \ \ \ \ \ \ \ \ \ \ \ } & 
	\mathrm{in} & (0,1).\text{ }%
	\end{array}%
	\right.  \label{stability}
	\end{equation}%
	The functions $u\in H_{\mathrm{loc}}^{1}(0,\infty ),$ and $f\in \mathcal{C}%
	([0,\infty ))$ in (\ref{control}) and (\ref{stability}) represent the
	control force and the feedback function respectively.
	
	\begin{figure}[H]
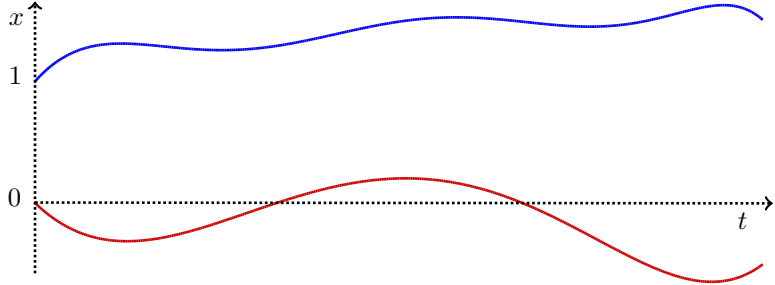

		\centering
		\definecolor{qqqqff}{rgb}{0,0,1} \definecolor{ccqqqq}{rgb}{0.8,0,0}  
		%
		\caption{The curve $(t,\protect\alpha(t))_{t\geq0}$ in red and $(t, \protect%
			\beta(t))_{t\geq0}$ in blue.}
	\end{figure}

	Controllability of system (\ref{control}) has been extensively studied in
	the recent past years; most of the papers dealt with the case of one moving
	endpoint with boundary conditions of the form 
	\begin{equation*}
	y(t,0)=0,\text{ \ }y(t,kt+1)=u(t),\text{ }k\in (0,1),\text{ }t\in (0,\infty
	).
	\end{equation*}
	
	In \cite{Cui1}, it has been shown that, with these boundary conditions,
	exact controllability holds for all times $T>\frac{e^{\frac{2k(k+1)}{1-k}}-1%
	}{2}.$ The same authors came back in \cite{Cui2} and improved the latter
	result to $T>\frac{e^{\frac{2k(k+1)}{\left( 1-k\right) ^{3}}}-1}{2}.$ Later,
	in \cite{Sun}, the controllability time has been improved to be $T>\frac{2}{%
		1-k}.$ In these papers, only a sufficient condition is provided for the
	exact controllability.
	
	Concerning the two moving endpoints case, the boundary functions considered
	in \cite{Sengouga} are of the form 
	\begin{equation*}
	\alpha (t)=-kt,\text{ \ }\beta (t)=rt+1,\text{ \ }t\in (0,\infty ),\text{ \ }%
	k,r\in \lbrack 0,1)\text{ with }r+k>0.
	\end{equation*}%
	It has been shown that exact controllability holds if, and only if $T\geq 
	\frac{2}{(1-k)(1-r)}$. More general boundary functions are considered in 
	\cite{Haak} with boundary conditions 
	\begin{equation*}
	y(t,0)=0,\text{ \ }y(t,s(t))=u(t),\text{ }t\in (0,\infty ),
	\end{equation*}%
	where $s:[0,\infty )\rightarrow (0,\infty )$ is assumed to be a $\mathcal{C}%
	^{1}$ function satisfying $\left\Vert s^{\prime }\right\Vert _{L^{\infty
		}(0,\infty )}<1.$ Furthermore, it has been assumed that $s$ must be in some
	admissible class of curves (see \cite{Haak} for more details). Under these
	assumptions, the authors proved that exact controllability holds if, and
	only if $T\geq s^{+}\circ \left( s^{-}\right) ^{-1}(0),$ where $s^{\pm
	}(t)=t\pm s(t)$. Also, they provided a controllability result when the
	control is located on the non-moving part of the boundary. By considering
	the boundary conditions 
	\begin{equation*}
	y(t,0)=u(t),\text{ \ }y(t,s(t))=0,\text{ }t\in (0,\infty ),
	\end{equation*}%
	they proved that exact controllability holds if, and only if $T\geq \left(
	s^{-}\right) ^{-1}(1).$ The same result has been proved in \cite{Gugat
		control} by using a different approach. In all the cited works, the proofs
	rely on the multipliers technique, the non-harmonic Fourier analysis or the
	d'Alembert solution of the wave equation.
	
	Recently, in \cite{Shao}, a new Carleman estimate has been established for
	the wave equation in non-cylindrical domains in more general settings. As a
	consequence, it has been shown for a boundary conditions as in (\ref{control}%
	) where $\alpha (t)<\beta (t),$ $t\in (0,\infty ),$ are smooth functions
	satisfying $\left\Vert \alpha ^{\prime }\right\Vert _{L^{\infty }(0,\infty
		)} $,$\left\Vert \beta ^{\prime }\right\Vert _{L^{\infty }(0,\infty )}<1,$
	that system (\ref{control}) is exactly controllable at time $T$ if $%
	T>T^{\ast }$ and not exactly controllable if $T<T^{\ast }$ where $T^{\ast }$
	is the required time by the geometric control condition, in other words, it
	is the time where a characteristic line with slope one emanating from the
	point $(0,0)$ hits the curve $\left( t,\beta (t)\right) _{t\geq 0}$ and
	reflected to intersect the curve $\left( t,\alpha (t)\right) _{t\geq 0}$ in
	the point $\left( T^{\ast },\alpha (T^{\ast })\right) $. Actually, this time
	can be computed explicitly in terms of the boundary curves, that is $T^{\ast
	}=$ $\left( \alpha ^{+}\right) ^{-1}\circ \beta ^{+}\circ \left( \beta
	^{-}\right) ^{-1}(0)$ where the functions $\alpha ^{\pm },\beta ^{\pm }$ are
	defined by $\alpha ^{\pm }(t)=t\pm \alpha (t),$ $\beta ^{\pm }(t)=t\pm \beta
	(t)$. However, the result doesn't cover the critical case $T=T^{\ast }.$
	
	As for the boundary stability of system (\ref{stability}) with nonautonomous
	damping, to the best of our knowlegde, the only existing result in the
	literature is in \cite{Ammari} where the authors dealt with the same system
	but with only one moving endpoint, i.e. 
	\begin{equation}
	y(t,0)=0,\text{ \ }y_{t}(t,a(t))+f(t)y_{x}(t,a(t))=0,\text{ }t\in (0,\infty
	),  \label{AAAA}
	\end{equation}%
	where $a$ is a strictly positive $1$-periodic function with $\left\Vert
	a^{\prime }\right\Vert _{L^{\infty }(0,\infty )}<1$ and $f$ is the feedback
	function. The authors proved exponential stablility of system (\ref%
	{stability}) for a particular class of feedbacks $f$. The proof relies on
	transforming problem (\ref{stability}) which is posed on non-cylindrical
	domain into a problem posed on cylindrical one, then making use of some
	known results of boundary stability of the 1D wave equation. If the damping
	function $f$ is constant and the boundary function $a$ is not periodic with
	derivative $\left\Vert a^{\prime }\right\Vert _{L^{\infty }(0,\infty )}<1,$
	it has been shown in \cite{Gugat stability} for $f=1$ that the solution vanishes at time $T$ for any $T\geq a^{+}\circ \left( a^{-}\right)
	^{-1}(0).$
	
	In this paper, we will improve all the previous results either for the
	boundary control or the boundary stability of the 1D wave equation by using
	the characteristics method. We shall build the unique exact solution to both
	systems (\ref{control}) and (\ref{stability}) in an appropriate energy
	space. To do so, we proceed by transforming both of systems to a first order
	hyperbolic system by introducing the Riemann invariants 
	\begin{equation}
	\left\{ 
	\begin{array}{cc}
	p= & y_{t}-y_{x,} \\ 
	q= & y_{t}+y_{x}.%
	\end{array}%
	\right.  \label{Riemann}
	\end{equation}%
	An elementary computation shows that system (\ref{control}) transforms into 
	\begin{equation}
	\left\{ 
	\begin{array}{ccc}
	p_{t}+p_{x}=0,\text{ \ \ \ \ \ \ \ \ \ \ \ \ \ \ \ \ \ \ \ \ \ \ \ \ \ \ \ \
		\ \ \ \ \ \ \ \ \ \ \ \ } & \mathrm{in} & Q,\text{ \ \ \ \ } \\ 
	q_{t}-q_{x}=0,\text{ \ \ \ \ \ \ \ \ \ \ \ \ \ \ \ \ \ \ \ \ \ \ \ \ \ \ \ \
		\ \ \ \ \ \ \ \ \ \ \ \ } & \mathrm{in} & Q,\text{ \ \ \ \ } \\ 
	\left( p+q\right) \left( t,\alpha (t)\right) =u^{\prime }(t),\text{ }\left(
	p+q\right) \left( t,\beta (t)\right) =0, & \mathrm{in} & (0,\infty ), \\ 
	p(0,x)=\widetilde{p}(x)\text{ },\text{ }q(0,x)=\widetilde{q}(x).\text{ \ \ \
		\ \ \ \ \ \ \ \ \ \ \ \ \ } & \mathrm{in} & (0,1).\text{ }%
	\end{array}%
	\right.  \label{control2}
	\end{equation}%
	In the same way, system (\ref{stability}) becomes 
	\begin{equation}
	\left\{ 
	\begin{array}{ccc}
	p_{t}+p_{x}=0,\text{ \ \ \ \ \ \ \ \ \ \ \ \ \ \ \ \ \ \ \ \ \ \ \ \ \ \ \ \
		\ \ \ \ \ \ \ \ \ \ \ \ \ } & \mathrm{in} & Q,\text{ \ \ \ \ } \\ 
	q_{t}-q_{x}=0,\text{ \ \ \ \ \ \ \ \ \ \ \ \ \ \ \ \ \ \ \ \ \ \ \ \ \ \ \ \
		\ \ \ \ \ \ \ \ \ \ \ \ \ } & \mathrm{in} & Q,\text{ \ \ \ \ } \\ 
	\left( p+F(t)q\right) \left( t,\alpha (t)\right) =0,\text{ }\left(
	p+q\right) \left( t,\beta (t)\right) =0, & \mathrm{in} & (0,\infty ), \\ 
	p(0,x)=\widetilde{p}(x)\text{ },\text{ }q(0,x)=\widetilde{q}(x),\text{ \ \ \
		\ \ \ \ \ \ \ \ \ \ \ \ \ } & \mathrm{in} & (0,1),\text{ }%
	\end{array}%
	\right.  \label{stability2}
	\end{equation}%
	where $F(t)=\frac{1-f(t)}{1+f(t)}$ with $1+f(t)\neq 0,$ $\forall t\geq 0.$
	
	Henceforth, we use the following notations:
	
	\begin{itemize}
		\item the spaces family $\left[ L^{2}(\alpha (t),\beta (t))\right] _{t\geq
			0} $ will be denoted by $L^{2}(\alpha (t),\beta (t)).$
		
		\item The spaces family $\left[ H_{\left( \beta (t)\right) }^{1}(\alpha
		(t),\beta (t))\right] _{t\geq 0}$ will be denoted by $H_{\left( \beta
			(t)\right) }^{1}(\alpha (t),\beta (t))$ where 
		\begin{equation*}
		H_{\left( \beta (t)\right) }^{1}(\alpha (t),\beta (t))=\left\{ h\in
		H^{1}(\alpha (t),\beta (t)),\text{ }h(\beta (t))=0,\text{ }t\geq 0\right\} .
		\end{equation*}
		
		\item For any function $z,$ the functions $z^{\pm }$ will represent the
		quantities $z^{\pm }(t)=t\pm z(t).$
		
		\item $C$ denotes a generic positive constant which might be different from
		line to line.
	\end{itemize}
	
	The Riemann coordinates introduced in (\ref{Riemann}) guarantee the
	equivalence of the transformed systems (\ref{control2}),(\ref{stability2}),
	with the original systems (\ref{control}),(\ref{stability}) up to an
	additive constant. All the results for the transformed systems will be
	proved in $\left[ L^{2}(\alpha (t),\beta (t))\right] ^{2},$ then the results
	for the original ones can be deduced by inverting the transformation.
	
	Since our approach consists in constructing the unique exact solutions to
	systems (\ref{control2}) and (\ref{stability2}), instead of studying each
	system separatly, we consider the following system 
	\begin{equation}
	\left\{ 
	\begin{array}{ccc}
	p_{t}+p_{x}=0,\text{ \ \ \ \ \ \ \ \ \ \ \ \ \ \ \ \ \ \ \ \ \ \ \ \ \ \ \ \
		\ \ \ \ \ \ \ \ \ \ \ \ \ \ \ \ \ } & \mathrm{in} & Q,\text{ \ \ \ \ } \\ 
	q_{t}-q_{x}=0,\text{ \ \ \ \ \ \ \ \ \ \ \ \ \ \ \ \ \ \ \ \ \ \ \ \ \ \ \ \
		\ \ \ \ \ \ \ \ \ \ \ \ \ \ \ \ \ } & \mathrm{in} & Q,\text{ \ \ \ \ } \\ 
	\left( p+F(t)q\right) \left( t,\alpha (t)\right) =v(t),\text{ }\left(
	p+q\right) \left( t,\beta (t)\right) =0, & \mathrm{in} & (0,\infty ), \\ 
	p(0,x)=\widetilde{p}(x)\text{ },\text{ }q(0,x)=\widetilde{q}(x).\text{ \ \ \
		\ \ \ \ \ \ \ \ \ \ \ \ \ \ \ \ \ \ } & \mathrm{in} & (0,1),\text{ }%
	\end{array}%
	\right.  \label{both}
	\end{equation}%
	where $v\in L_{\mathrm{loc}}^{2}(0,\infty )$ stands for $u^{\prime }.$ Note
	that if $F\equiv 1$ then system (\ref{both}) turns to be (\ref{control2}),
	and if $v\equiv 0,$ system (\ref{both}) turns to be (\ref{stability2}).
	Observe that the solutions to the first and the second equations of (\ref%
	{both}) satisfy 
	\begin{equation}
	\frac{d}{dt}p(t,c+t)=\frac{d}{dt}q(t,c-t)=0,\text{ \ }t\geq 0,\text{ }c\in 
	\mathbb{R}
	,  \label{charac}
	\end{equation}%
	hence, $p$ (res. $q$) is constant along the characteristic lines $x-t=c$
	(resp. $x+t=c$). The idea is to use the boundary conditions 
	\begin{equation}
	\left( p+F(t)q\right) \left( t,\alpha (t)\right) =v(t),\text{ }\left(
	p+q\right) \left( t,\beta (t)\right) =0,\text{ }t>0,  \label{boundary}
	\end{equation}%
	and the reflection of the characteristic lines $x\pm t=c$, $c\in 
	\mathbb{R}
	,$ on the boundary curves $\left( t,\alpha (t)\right) _{t\geq 0}$ and $%
	\left( t,\beta (t)\right) _{t\geq 0}$ to find the unique solution to system (%
	\ref{both}). Along this work, we assume that the boundary curves satisfy 
	\begin{equation}
	\text{ }\alpha (t)<\beta (t),\text{ }%
	\forall t>0,\text{ }\alpha ,\beta \in \mathcal{C}^{1}(0,\infty )\text{ },%
	\text{ }\left\Vert \alpha ^{\prime }\right\Vert _{L^{\infty }(%
		\mathbb{R}
		_{+})},\left\Vert \beta ^{\prime }\right\Vert _{L^{\infty }(%
		\mathbb{R}
		_{+})}<1.  \label{assumption}
	\end{equation}%
	The size assumption in (\ref{assumption}) guarantees that the characteristic
	lines $x=t+c$ (resp. $x=c-t$) meet the curve $\left( t,\alpha (t)\right)
	_{t\geq 0}$ (resp. $\left( t,\beta (t)\right) _{t\geq 0})$ in finite time;
	also, they serve to ensure that the characteristic lines $x\pm t=c$ are not
	gliding on the boundary curves or are not out of $Q.$ In fact, assumption (%
	\ref{assumption}) is necessary for the existence of solutions. A
	straightforward consequence of assumption (\ref{assumption}) is that the
	functions $\alpha ^{\pm }:[0,\infty )\rightarrow \lbrack 0,\infty )$ and $%
	\beta ^{\pm }:[0,\infty )\rightarrow \lbrack \pm 1,\infty )$ are invertible.
	In the sequel, we use the standard notations to denote their inverses by $%
	\left( \alpha ^{\pm }\right) ^{-1}$ and $\left( \beta ^{\pm }\right) ^{-1}.$
	\begin{figure}[H]
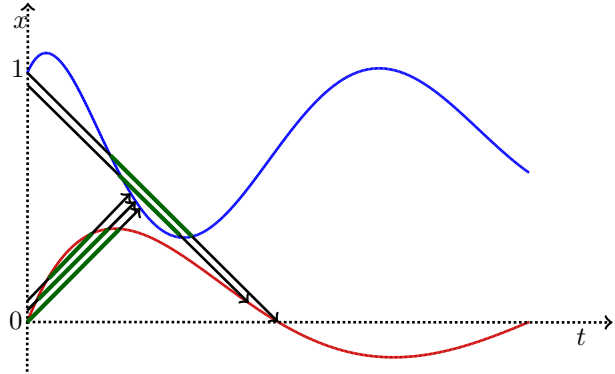

		\centering
		\definecolor{qqwuqq}{rgb}{0,0.39215686274509803,0}  %
		\definecolor{ccqqqq}{rgb}{0.8,0,0} \definecolor{qqqqff}{rgb}{0,0,1}  

		\caption{An example of a boundary curves $(t,\protect\alpha(t)))_{t\geq0}$
			and $(t,\protect\beta(t)))_{t\geq0}$ that do not satisfy assumption (\protect
			\ref{assumption}). The values of the solution are not defined on the green
			part of the characteristic lines lying under or above these curves.}
	\end{figure}

	\section{Main results}
	
	We start by giving the well-posedness result for system (\ref{both}).
	
	\begin{theorem}
		Let $\left( \widetilde{p},\widetilde{q},v,F\right) \in \left[ L^{2}(0,1)%
		\right] ^{2}\times L_{\mathrm{loc}}^{2}(0,\infty )\times \mathcal{C}%
		([0,\infty ))$. Assume that the boundary curves $(t,\alpha (t))_{t\geq 0}$
		and $(t,\beta (t))_{t\geq 0}$ satisfy (\ref{assumption}). Then, there exists
		a unique solution to system (\ref{both}) satisfying 
		\begin{equation}
		(p,q)\in \mathcal{C}\left( 0,t;\left[ L^{2}(\alpha (t),\beta (t))\right]
		^{2}\right) ,\text{ }t\geq 0.  \label{regularity}
		\end{equation}
	\end{theorem}
	
	The proof of this theorem is a straightforward consequence of the explicit
	construction of the unique solution that will be done in Section \ref%
	{section construction}.
	
	\begin{remark}
		By inverting the transformation given in (\ref{Riemann}), we obtain 
		\begin{equation*}
		y_{t}=\frac{p+q}{2}\text{ \ , \ }y_{x}=\frac{q-p}{2},
		\end{equation*}%
		hence, for any $\left( y_{0},y_{1},u,f\right) \in H_{\left( 1\right)
		}^{1}(0,1)\times L^{2}(0,1)\times H_{loc}^{1}(0,\infty )\times \mathcal{C}%
		([0,\infty )),$ the solutions to systems (\ref{control}) and (\ref{stability}%
		) satisfy the regularity 
		\begin{equation*}
		y\in C\left( 0,t;H_{\left( \beta (t)\right) }^{1}(\alpha (t),\beta
		(t)\right) \cap C^{1}\left( 0,t;L^{2}(\alpha (t),\beta (t)\right) ,\text{ }%
		t\geq 0.
		\end{equation*}
	\end{remark}
	
	\subsection{Controllability result}
	
	\begin{definition}
		System (\ref{control}) is said to be exactly controllable at time $T>0$ if
		for any initial state $(y_{0},y_{1})\in H_{\left( 1\right) }^{1}(0,1)\times
		L^{2}(0,1)$ and for any target state $(h,k)\in H_{\left( \beta (T)\right)
		}^{1}(\alpha (T),\beta (T))\times L^{2}(\alpha (T),\beta (T)),$ there exists
		a control $u\in H_{\mathrm{loc}}^{1}(0,\infty )$ such that $\left(
		y(T),y_{t}(T)\right) =(h,k).$
	\end{definition}
	
	The following result shows that the minimal time $T^{\ast }$ where exact
	controllability is possible depends on the movement of the boundaries and
	can be represented explicitly in terms of the functions $\alpha ^{\pm }$ and 
	$\beta ^{\pm }$. Moreover, also the uniqe exact control for $T^{\ast }$ can
	be represented explictily using these functions.
	
	\begin{theorem}
		\label{theorem control}Let $(y_{0},y_{1})\in H_{\left( 1\right)
		}^{1}(0,1)\times L^{2}(0,1).$ Assume that the boundary curves $(t,\alpha
		(t))_{t\geq 0}$ and $(t,\beta (t))_{t\geq 0}$ satisfy (\ref{assumption}).
		System (\ref{control}) is exactly controllable at time $T>0$ if, and only if 
		$T\geq T^{\ast }=$ $\left( \alpha ^{+}\right) ^{-1}\circ \beta ^{+}\circ
		\left( \beta ^{-}\right) ^{-1}(0).$ Further, if $T=T^{\ast },$ there exists
		a unique control $u\in H^{1}(0,T^{\ast })$ steering the solution $(y,y_{t})$
		to system (\ref{control}) to the equilibrium point $(0,0)$ given by 
		\begin{equation}
		u(t)=\left\{ 
		\begin{array}{ccc}
		\int_{0}^{t}y_{1}\left( \alpha ^{+}(s)\right) ds+y_{0}\left( \alpha
		^{+}(t)\right) , & \mathrm{if} & t\in \left[ 0,\left( \alpha ^{+}\right)
		^{-1}(1)\right) ,\text{ } \\ 
		\begin{array}{c}
		\\ 
		y_{0}\left( -\beta ^{-}\circ \left( \beta ^{+}\right) ^{-1}\circ \alpha
		^{+}(t)\right) \\ 
		+\int_{0}^{\left( \alpha ^{+}\right) ^{-1}(1)}y_{1}\left( \alpha
		^{+}(s)\right) ds \\ 
		-\int_{\left( \alpha ^{+}\right) ^{-1}(1)}^{t}y_{1}\left( -\beta ^{-}\circ
		\left( \beta ^{+}\right) ^{-1}\circ \alpha ^{+}(s)\right) ds,%
		\end{array}
		& \mathrm{if} & t\in \left[ \left( \alpha ^{+}\right) ^{-1}(1),T^{\ast
		}\right) ,%
		\end{array}%
		\right.  \label{co}
		\end{equation}
	\end{theorem}
	
	\begin{remark}
		The controllability result still makes sense even if the boundary curves $%
		(t,\alpha (t))_{t\geq 0}$ and $(t,\beta (t))_{t\geq 0}$ are allowed to
		intersect in time larger than $T^{\ast }.$
	\end{remark}
	
	\begin{remark}
		Let us consider the particular case $\alpha (t)=kt,$ $\beta (t)=rt+1$, $%
		k,r\in (-1,1),$ with $\frac{2(k-r)}{(1-r)(1+k)}<\frac{1}{2}$ (The last
		assumption guarantees that the boundary curves do not intersect before $%
		T^{\ast }$)$.$ In this case, it can be checked that $T^{\ast }$ is given by $%
		T^{\ast }=\frac{2}{(1-r)(k+1)}$ which is the same time found in \cite%
		{Sengouga}$.$ In particular, if $\alpha \equiv 0$ and $\beta \equiv 1$, we
		obtain the classical result $T^{\ast }=2.$
	\end{remark}
	
	\begin{remark}
		The minimal time $T^{\ast }$ is precisely the necessary time for the main
		characteristic line issued from the point $(0,0)$ to touch again the curve $%
		\left( t,\alpha (t)\right) _{t\geq 0}$ in the point $\left( T^{\ast },\alpha
		(T^{\ast })\right) $ after having been reflected from the curve $\left(
		t,\beta (t)\right) _{t\geq 0}$. More precisely, the characteristic line $x=t$
		hits the curve $\left( t,\beta (t)\right) _{t\geq 0}$ in the point $\left(
		\left( \beta ^{-}\right) ^{-1}(0),\beta \left( \left( \beta ^{-}\right)
		^{-1}(0)\right) \right) $. The reflected characteristic line passing through
		the last point, i.e. $x=-t+\beta ^{+}\circ \left( \beta ^{-}\right) ^{-1}(0)$
		hits the curve $\left( t,\alpha (t)\right) _{t\geq 0}$ in the point $\left(
		T^{\ast },\alpha (T^{\ast })\right) $. If the control $u$ is located on the
		curve $\left( t,\beta (t)\right) _{t\geq 0}$ instead of $\left( t,\alpha
		(t)\right) _{t\geq 0}$, then $T^{\ast \ast }$ is the analogous time for the
		main characteirstic line issued form the point $(0,1)$ with negative slope$.$
		In this case $T^{\ast \ast }=\left( \beta ^{-}\right) ^{-1}\circ \alpha
		^{-}\circ \left( \alpha ^{+}\right) ^{-1}(1)$.
	\end{remark}
	
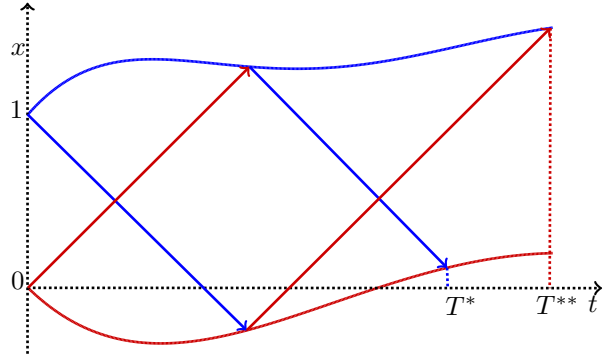
\begin{figure}[H]
		\centering
		\definecolor{qqqqff}{rgb}{0,0,1} \definecolor{ccqqqq}{rgb}{0.8,0,0}  

	\end{figure}

	\subsection{Stability result}
	
	For the sake of lighting notations, we introduce the function $\phi :=\phi
	(\alpha ,\beta )$ defined by 
	\begin{equation}
	\phi :=\alpha ^{-}\circ \left( \alpha ^{+}\right) ^{-1}\circ \beta ^{+}\circ
	\left( \beta ^{-}\right) ^{-1}.  \label{phi}
	\end{equation}
	By assumption (\ref{assumption}), the function $\phi :[-1,\infty
	)\rightarrow \lbrack \alpha ^{-}\circ \left( \alpha ^{+}\right)
	^{-1}(1),\infty )$ is well defined and increasing function as composition of
	increasing functions, and hence invertible with inverse 
	\begin{equation*}
	\phi ^{-1}:=\beta ^{-}\circ \left( \beta ^{+}\right) ^{-1}\circ \alpha
	^{+}\circ \left( \alpha ^{-}\right) ^{-1}.
	\end{equation*}
	
	Let $\left( \psi _{n}\right) _{n\geq 0}$ be a sequence of functions such
	that 
	\begin{eqnarray}
	\psi _{n} &:&[0,\phi (0))\rightarrow \lbrack 0,\infty )  \label{psi} \\
	\tau &\mapsto &\psi _{n}(\tau )=\prod\limits_{i=0}^{n}\left\vert F\left(
	\left( \alpha ^{-}\right) ^{-1}\circ \phi ^{\left[ i\right] }(\tau )\right)
	\right\vert .  \notag
	\end{eqnarray}
	The notation $\phi ^{\left[ n\right] }$ refeers to the $n^{th}$ composed of $%
	\phi ^{\left[ n\right] }$ i.e.
	
	\begin{equation*}
	\phi ^{\left[ n\right] }=\underset{n\text{ times}}{\underbrace{\phi \circ
			\phi \circ \phi \circ \phi \circ \cdot \cdot \cdot \circ \phi }},
	\end{equation*}%
	with the convention $\phi ^{\left[ 0\right] }=I.$ The following result shows
	that the asymptotic behaviour of the solution to system (\ref{stability})
	for large time relies heavily on the behaviour of the sequence of functions $%
	\left( \psi _{n}(\tau )\right) _{n\geq 0}$ defined in (\ref{psi}) when $%
	n\longrightarrow \infty .$
	
	\begin{theorem}
		\label{theorem stab}Let $(y_{0},y_{1})\in H_{\left( 1\right)
		}^{1}(0,1)\times L^{2}(0,1).$ Assume that the boundary curves $(t,\alpha
		(t))_{t\geq 0}$ and $(t,\beta (t))_{t\geq 0}$ satisfy (\ref{assumption}). In
		addition, assume that 
		\begin{equation}
		\phi (\tau )<\cdot \cdot \cdot <\phi ^{\left[ n\right] }(\tau )<\phi ^{\left[
			n+1\right] }(\tau )\underset{n\rightarrow \infty }{\longrightarrow }\infty ,%
		\text{ }\forall \tau \in \lbrack 0,\phi (0)),  \label{increasing}
		\end{equation}%
		then,%
		\begin{equation*}
		\left\Vert \left( y(t),y_{t}(t)\right) \right\Vert _{H_{\left( \beta
				(t)\right) }^{1}(\alpha (t),\beta (t))\times L^{2}(\alpha (t),\beta (t))}%
		\underset{t\rightarrow \infty }{\longrightarrow }0,
		\end{equation*}%
		if, and only if 
		\begin{equation}
		\psi _{n}(\tau )\underset{n\rightarrow \infty }{\longrightarrow }0,\text{ }%
		\forall \tau \in \lbrack 0,\phi (0)).  \label{assum}
		\end{equation}%
		If there exists $g\in \mathcal{C}(%
		\mathbb{R}
		,(0,\infty ))$\ such that 
		\begin{equation}
		\psi _{n}(\tau )\underset{n\rightarrow \infty }{\sim }Cg\left( \phi ^{\left[
			n\right] }(\tau )\right) ,\text{ }\forall \tau \in \lbrack 0,\phi (0)),
		\label{g}
		\end{equation}%
		then, the solution to system (\ref{stability}) decays like $g(t),$ i.e.
		there exists a positive constant $C$ such that 
		\begin{equation}
		\left\Vert \left( y(t),y_{t}(t)\right) \right\Vert _{H_{\left( \beta
				(t)\right) }^{1}(\alpha (t),\beta (t))\times L^{2}(\alpha (t),\beta
			(t))}\leq Cg(t)\left\Vert \left( y_{0},y_{1}\right) \right\Vert _{H_{\left(
				1\right) }^{1}(0,1)\times L^{2}(0,1)}.  \label{rate}
		\end{equation}%
		In particular, the solution to system (\ref{stability}) $\left(
		y(t),y_{t}(t)\right) $ decays exponentially to zero with growth $\omega >0$,
		i.e. there exists $M\geq 1$ such that 
		\begin{equation*}
		\left\Vert \left( y(t),y_{t}(t)\right) \right\Vert _{H_{\left( \beta
				(t)\right) }^{1}(\alpha (t),\beta (t))\times L^{2}(\alpha (t),\beta
			(t))}\leq Me^{-t\omega }\left\Vert \left( y_{0},y_{1}\right) \right\Vert
		_{H_{\left( 1\right) }^{1}(0,1)\times L^{2}(0,1)},\text{ }\forall t\geq 0,
		\end{equation*}%
		if, and only if 
		\begin{equation}
		\underset{\tau \in \lbrack 0,\phi (0))}{\sup }\underset{n\rightarrow \infty }%
		{\lim }\frac{\ln \psi _{n}(\tau )}{\phi ^{\left[ n\right] }(\tau )}=-\omega .
		\label{A1}
		\end{equation}%
		If $f\equiv 1$, the solution to system (\ref{stability}) vanishes in finite
		time $T$ if, and only if $T\geq T^{\ast }=$ $\left( \alpha ^{+}\right)
		^{-1}\circ \beta ^{+}\circ \left( \beta ^{-}\right) ^{-1}(0),$ i.e. 
		\begin{equation*}
		y(T)\equiv y_{t}(T)\equiv 0,\text{ }\forall T\geq T^{\ast }=\left( \alpha
		^{+}\right) ^{-1}\circ \beta ^{+}\circ \left( \beta ^{-}\right) ^{-1}(0).
		\end{equation*}
	\end{theorem}
	
	Let us illustrate the previous theorem by some examples.
	
	\begin{example}[Cylindrical domain]
		\label{example}If $Q$ is cylindrical domain, i.e. $\alpha \equiv 0$ and $%
		\beta \equiv 1,$ the function $\phi $ defined in (\ref{phi}) is given by $%
		\phi (\tau )=\tau +2,$ then, $\phi ^{\left[ n\right] }(\tau )=\tau +2n.$
		Therefore, the functions sequence $\left( \psi _{n}\right) _{n\geq 0}$
		defined in (\ref{psi}) takes the form 
		\begin{eqnarray}
		\psi _{n} &:&[0,2)\rightarrow \lbrack 0,\infty )  \label{ksii} \\
		\tau &\mapsto &\psi _{n}(\tau )=\prod\limits_{i=0}^{n}\left\vert F\left(
		\tau +2i\right) \right\vert ,  \notag
		\end{eqnarray}%
		In this case, assumptions of Theorem (\ref{theorem stab}) can be checked
		easily. Note that since system (\ref{stability}) is non-autonomous ($f$ is
		time dependent), the decay rate is not necessarily exponential. Below, we
		illustarte this fact by several examples:
		
		\begin{itemize}
			\item Exponential decay:
			
			Let $f(t)=\frac{2-\sin (\pi t)}{2+\sin (\pi t)},$ therefore, $F(t)=\frac{
				\sin (\pi t)}{2},$ thus, 
			\begin{equation*}
			\psi _{n}(\tau )=\prod\limits_{i=0}^{n}\left\vert F\left( \tau +2i\right)
			\right\vert =\left[ \frac{\sin (\pi \tau )}{2}\right] ^{n+1}.
			\end{equation*}
			By (\ref{A1}), we have 
			\begin{eqnarray*}
				\underset{\tau \in (0,1)\cup (1,2)}{\sup }\underset{n\rightarrow \infty }{
					\lim }\frac{\ln \psi _{n}(\tau )}{\phi ^{\left[ n\right] }(\tau )} &=& 
				\underset{\tau \in (0,1)\cup (1,2)}{\sup }\underset{n\rightarrow \infty }{
					\lim }\frac{\left( n+1\right) \ln \left\vert \frac{\sin (\pi \tau )}{2}
					\right\vert }{\tau +2n} \\
				&=&\underset{\tau \in (0,1)\cup (1,2)}{\sup }\frac{1}{2}\ln \left\vert \frac{
					\sin (\pi \tau )}{2}\right\vert =-\frac{\ln 2}{2},
			\end{eqnarray*}
			therefore, exponential decay occurs with growth bound $\omega =\frac{\ln 2}{
				2 }.$
			
			\item Polynomial decay:
			
			Let $f(t)=\frac{\left( t+1\right) ^{-s}-(t+3)^{-s}}{\left( t+1\right)
				^{-s}+(t+3)^{-s}},$ $s>0,$ then $F(t)=\left( \frac{t+3}{t+1}\right) ^{-s},$
			consequently, the functions sequence $\left( \psi _{n}\right) _{n\geq 0}$
			defined in (\ref{ksii}) takes the form 
			\begin{equation*}
			\psi _{n}(\tau )=\prod\limits_{i=0}^{n}\left\vert F\left( \tau +2i\right)
			\right\vert =\prod\limits_{i=0}^{n}\left\vert \left( \frac{\tau +2i+3}{\tau
				+2i+1}\right) ^{-s}\right\vert =\left( \frac{\tau +2n+3}{\tau +1}\right)
			^{-s}.
			\end{equation*}
			Set $g(t)=\left( t+1\right) ^{-s},$ $s>0.$ A simple computation shows that 
			\begin{equation*}
			\underset{n\rightarrow \infty }{\lim }\frac{\psi _{n}(\tau )}{g(\phi
				^{\lbrack n]})}=\underset{n\rightarrow \infty }{\lim }\frac{\psi _{n}(\tau ) 
			}{g(\tau +2n)}=\frac{1}{\tau +1},\text{ }\tau \in \lbrack 0,2),\text{ }
			\end{equation*}
			thus, by (\ref{rate}), the solution to system (\ref{stability}) decays like $%
			(t+1)^{-s},$ $s>0.$
			
			\item Logarithmic decay:
			
			Let $f(t)=\frac{\log ^{-s}(t+1)-\log ^{-s}(t+3)}{\log ^{-s}(t+1)+\log
				^{-s}(t+3)},$ $s>0,$ then $F(t)=\left( \frac{\log (t+3)}{\log (t+1)}\right)
			^{-s},$ consequently, we obtain 
			\begin{equation*}
			\psi _{n}(\tau )=\prod\limits_{i=0}^{n}\left\vert F\left( \tau +2i\right)
			\right\vert =\prod\limits_{i=0}^{n}\left\vert \left( \frac{\log (\tau +2i+3) 
			}{\log (\tau +2i+1)}\right) ^{-s}\right\vert =\left\vert \left( \frac{\log
				(\tau +2n+3)}{\log (\tau +1)}\right) ^{-s}\right\vert .
			\end{equation*}
			\ By letting $g(t)=\log ^{-s}(t+1),$ $s>0,$ we get 
			\begin{equation*}
			\underset{n\rightarrow \infty }{\lim }\frac{\psi _{n}(\tau )}{g(\phi
				^{\lbrack n]})}=\log ^{s}(\tau +1)\underset{n\rightarrow \infty }{\lim }
			\left\vert \left( \frac{\log (\tau +2n+3)}{\log (\tau +2n+1)}\right)
			^{-s}\right\vert =\log ^{s}(\tau +1),\text{ }\tau \in \lbrack 0,2),
			\end{equation*}
			hence, (\ref{rate}) is satisfied with $g(t)=\log ^{-s}(t+1),$ $s>0.$
			
			\item Stability in finite time:
			
			Let $f(t)=\frac{t}{2+t},$ therefore, $F(t)=\frac{1}{t+1},$ consequently, we
			obtain 
			\begin{equation*}
			\psi _{n}(\tau )=\prod\limits_{i=0}^{n}\frac{1}{\tau +2i+1}=\frac{1}{\left(
				\tau +1\right) 2^{n}n!}\prod\limits_{i=1}^{n}\left( \frac{\tau +1}{2i}%
			+1\right) ^{-1}.
			\end{equation*}
			
			A simple computation shows that 
			\begin{equation*}
			\log \prod\limits_{i=1}^{n}\left( \frac{\tau +1}{2i}+1\right) ^{-1}\underset{%
				n\rightarrow \infty }{\sim }C(\tau )\log n^{-\frac{\tau +1}{2}},
			\end{equation*}
			
			where $C(\tau )$ is a positive constant depending on $\tau .$ So, we get 
			\begin{equation}
			\psi _{n}(\tau )\underset{n\rightarrow \infty }{\sim }\frac{C(\tau )}{\left(
				\tau +1\right) n^{\frac{\tau +1}{2}}2^{n}n!},  \label{poch}
			\end{equation}%
			which by (\ref{assum}) implies that the solution to system (\ref{stability})
			decays to zero. To prove stability in finite time it suffices to show that
			the gowth bound defined in (\ref{A1}) is infinite. By using (\ref{poch}) we
			obtain 
			\begin{equation*}
			\underset{n\rightarrow \infty }{\lim }\frac{\ln \psi _{n}(\tau )}{\phi ^{%
					\left[ n\right] }(\tau )}=-\underset{n\rightarrow \infty }{\lim }\frac{\ln
				n^{\frac{\tau +1}{2}}+n\ln 2+\ln n!}{2n+\tau }=-\infty =-\omega .
			\end{equation*}%
			This phenomena is due to the fact that $f(t)\underset{t\rightarrow \infty }{%
				\longrightarrow }1$ which by Theorem \ref{theorem stab} leads to stability
			in finite time.
		\end{itemize}
	\end{example}
	
	\begin{example}[Non cylindrical domain]
		\label{example2}Things are more delicate in the non-cylindrical case.
		Consider a boundary functions of the form $\alpha (t)=rt,$ $\beta (t)=kt+1,$ 
		$r,k\in (-1,1).$ To guarantee that $\alpha (t)\neq \beta (t),\forall t\geq 0,$
		we assume that $k\geq r.$ The function $\phi $ defined in (\ref{phi}) will
		be given by 
		\begin{equation*}
		\phi (\tau )=\frac{\left( 1+k\right) \left( 1-r\right) }{\left( 1-k\right)
			\left( 1+r\right) }\tau +\frac{2\left( 1-r\right) }{\left( 1-k\right) \left(
			1+r\right) }=a\tau +b,
		\end{equation*}
		therefore, we obtain 
		\begin{equation*}
		\phi ^{\lbrack n]}(\tau )=\left\{ 
		\begin{array}{ccc}
		a^{n}\left( \tau -\frac{b}{1-a}\right) +\frac{b}{1-a}, & \mathrm{if} & r<k,
		\\ 
		\begin{array}{c}
		\\ 
		\tau +\frac{2n}{1+r},\text{ \ \ \ \ \ \ \ \ \ \ \ \ \ }%
		\end{array}
		& 
		\begin{array}{c}
		\\ 
		\mathrm{if}%
		\end{array}
		& 
		\begin{array}{c}
		\\ 
		r=k.%
		\end{array}%
		\end{array}
		\right.
		\end{equation*}
		Consequently, 
		\begin{equation}
		\left( \alpha ^{-}\right) ^{-1}\circ \phi ^{\lbrack n]}(\tau )=\left\{ 
		\begin{array}{ccc}
		a^{n}\left( \frac{\tau }{1-r}-\frac{b}{\left( 1-a\right) \left( 1-r\right) }
		\right) +\frac{b}{\left( 1-a\right) \left( 1-r\right) }, & \mathrm{if} & r<k,
		\\ 
		\begin{array}{c}
		\\ 
		\frac{\tau }{1-r}+\frac{2n}{\left( 1+r\right) \left( 1-r\right) },\text{ \ \
			\ \ \ \ \ \ \ \ \ \ \ \ \ \ \ \ \ \ \ \ }%
		\end{array}
		& 
		\begin{array}{c}
		\\ 
		\mathrm{if}%
		\end{array}
		& 
		\begin{array}{c}
		\\ 
		r=k.%
		\end{array}%
		\end{array}
		\right.  \label{AA}
		\end{equation}
		For simplicity, let us take $f$ as in the previous example, $f(t)=\frac{t}{
			2+t}$ which implies that $F(t)=\frac{1}{t+1}.$ So, we have:
		
		\begin{itemize}
			\item If $r<k:$
			
			From (\ref{AA}), we can check that (\ref{increasing}) is satisfied if, and
			only if $a>1.$ To verify (\ref{assum}), it is enough to estimate its
			asymptotics. So, we have 
			\begin{eqnarray*}
				\psi _{n}(\tau ) &=&\prod\limits_{i=0}^{n}\frac{1}{\left\vert a^{i}\left( 
					\frac{\tau }{1-r}-\frac{b}{\left( 1-a\right) \left( 1-r\right) }\right) +%
					\frac{b}{\left( 1-a\right) \left( 1-r\right) }+1\right\vert }%
				=\prod\limits_{i=0}^{n}\frac{1}{\left\vert a^{i}s(\tau )+z\right\vert } \\
				&=&\frac{1}{a^{\frac{n(n+1)}{2}}s^{n+1}(\tau )}\prod\limits_{i=0}^{n}\left%
				\vert 1+\frac{z}{a^{i}s(\tau )}\right\vert ^{-1}.
			\end{eqnarray*}%
			Since $a>1,$ the series $\sum_{i=0}^{\infty }\ln \left( 1+\frac{z}{%
				a^{i}s(\tau )}\right) $ converges, we obtain, 
			\begin{equation}
			\psi _{n}(\tau )\underset{n\rightarrow \infty }{\sim }C(r,k,\tau )a^{-\frac{%
					n(n+1)}{2}}s^{-n-1}(\tau ),\text{ }\forall \tau \in \left[ 0,b\right) ,
			\label{pro}
			\end{equation}%
			where $C(r,k,\tau )$ is a positive constant depending on $r,k$ and $\tau .$
			In view of (\ref{pro}), assumption (\ref{assum}) is satisfied if $a>1,$
			namely, 
			\begin{equation}
			\frac{\left( 1+k\right) \left( 1-r\right) }{\left( 1-k\right) \left(
				1+r\right) }>1,  \label{k}
			\end{equation}%
			i.e. the solution to system (\ref{stability}) decays to zero if (\ref{k}) is
			satisfied. On the contrary of the cylindrical domain case, even with this
			choice of the feedback function $f,$ (\ref{A1}) is not satisfied, and hence,
			exponential stability cannot occur. Indeed, we have 
			\begin{equation*}
			\underset{n\rightarrow \infty }{\lim }\frac{\ln \psi _{n}(\tau )}{\phi ^{%
					\left[ n\right] }(\tau )}=-\underset{n\rightarrow \infty }{\lim }\frac{n^{2}%
			}{2a^{n}}\ln \left\vert a\right\vert =0,\text{ }\forall \tau \in \left[
			0,b\right) .
			\end{equation*}
			
			We still be able to get an idea about the decay rate. From (\ref{pro}), we
			observe that the term that really matters is $a^{-\frac{n^{2}}{2}},$ so, for 
			$g(t)=e^{-\frac{1}{2}\log _{a}^{2}(t)}$, we obtain 
			\begin{equation*}
			a^{-\frac{n^{2}}{2}}\underset{n\rightarrow \infty }{\sim }Cg(a^{n}s(\tau
			)+z),\text{ }\forall \tau \in \left[ 0,b\right) .
			\end{equation*}%
			Note that we did not lose too much since $g$ decays to zero faster than any
			polynomial function. This loss can be justified by the fact that the
			characteristic lines will need a larger time to reflect on the two endpoints
			when $t$ becomes larger.
			
			\item If $k=r:$
			
			In this case, the lines $x=rt$ and $x=kt+1$ are parallel, therefore, the
			characteristic speeds are the same for all time, so we might expect
			stability in finite time with this choice of $f$. Let us first check that whether
			the solution to system (\ref{stability}) decays exponentially or not. By (\ref{poch}%
			), the functions sequence $\left( \psi _{n}\right) _{n\geq 1}$ behaves like 
			\begin{eqnarray*}
				\psi _{n}(\tau ) &=&\prod\limits_{i=0}^{n}\left\vert \frac{1}{\frac{\tau }{%
						1-r}+1+\frac{2i}{\left( 1+r\right) \left( 1-r\right) }}\right\vert \\
				&=&\frac{\left( 1+r\right) ^{n}\left( 1-r\right) ^{n}}{2^{n}n!\left( \frac{%
						\tau }{1-r}+1\right) }\prod\limits_{i=1}^{n}\left\vert \frac{1}{\frac{\left(
						1+r\right) \left( 1-r\right) }{2i}\left( \frac{\tau }{1-r}+1\right) +1}%
				\right\vert . \\
				&&\underset{n\rightarrow \infty }{\sim }C(r,\tau )\frac{\left( 1+r\right)
					^{n}\left( 1-r\right) ^{n+1}}{2^{n}n!\left( \tau +1-r\right) n^{\frac{\left(
							1+r\right) \left( \tau +1-r\right) }{2}}},
			\end{eqnarray*}%
			where $C(r,\tau )$ is a positive constant depending on $r$ and $\tau .$ By
			using (\ref{A1}), we get 
			\begin{equation*}
			\underset{n\rightarrow \infty }{\lim }\frac{\ln \psi _{n}(\tau )}{\phi ^{%
					\left[ n\right] }(\tau )}=\underset{n\rightarrow \infty }{\lim }-\frac{\ln n!%
			}{2n}=-\infty =-\omega ,
			\end{equation*}%
			therefore, the solution to system (\ref{stability}) will vanish in finite
			time. This is due to the fact that $f(t)\underset{t\rightarrow \infty }{%
				\longrightarrow }1.$
		\end{itemize}
	\end{example}
	
	\begin{example}[Constant feedback]
		\label{example3}Consider the case when $f$ is a constant such that $f\neq 1$
		with keeping $\alpha $ and $\beta $ as in the previous example. A simple
		computation yields 
		\begin{equation*}
		\psi _{n}=F^{n+1}=\left\vert \frac{f-1}{f+1}\right\vert ^{n+1}.
		\end{equation*}%
		Therefore, by using the formula in (\ref{A1}), we arrive at:
		
		\begin{itemize}
			\item If $r<k:$
			
			We can check that the decay is not exponential. Indeed, 
			\begin{equation*}
			\underset{n\rightarrow \infty }{\lim }\frac{\ln \psi _{n}(\tau )}{\phi ^{%
					\left[ n\right] }(\tau )}=\underset{n\rightarrow \infty }{\lim }\frac{\left(
				n+1\right) \ln \left\vert \frac{f-1}{f+1}\right\vert }{a^{n}}=0,\text{ }%
			\forall \tau \in \left[ 0,b\right) .
			\end{equation*}%
			Nonetheless, by (\ref{rate}), we can determine the decay rate for a
			particular values of $f$ . Let $g(t)=t^{-s}.$ It is easy to check that if $%
			a^{-s}=\frac{f-1}{f+1}$ for some $s>0$ then 
			\begin{equation*}
			\underset{n\rightarrow \infty }{\lim }\frac{\psi _{n}(\tau )}{g\left( \phi ^{%
					\left[ n\right] }(\tau )\right) }=\underset{n\rightarrow \infty }{\lim }%
			\frac{\left\vert \frac{f-1}{f+1}\right\vert ^{n+1}}{\left( a^{n}s(\tau
				)+z\right) ^{-s}}=C(\tau ,r,k),\text{ }\forall \tau \in \left[ 0,b\right) ,
			\end{equation*}%
			where $C(r,k,\tau )$ is a positive constant depending on $r,k$ and $\tau .$
			Hence, the solution decays like $t^{-s},$ $s>0.$
			
			\item If $r=k:$
			
			In this case, we have 
			\begin{equation*}
			\underset{n\rightarrow \infty }{\lim }\frac{\ln \psi _{n}(\tau )}{\phi ^{%
					\left[ n\right] }(\tau )}=\underset{n\rightarrow \infty }{\lim }\frac{\left(
				n+1\right) \ln \left\vert \frac{f-1}{f+1}\right\vert }{\frac{2n}{\left(
					1+r\right) \left( 1-r\right) }}=\frac{\left( 1+r\right) \left( 1-r\right) }{2%
			}\ln \left\vert \frac{f-1}{f+1}\right\vert =-\omega ,
			\end{equation*}%
			hence, exponential decay occurs with growth bound $\omega .$ In particular,
			if $Q$ is a cylindrical domain $(r=0),$ system (\ref{stability}) is
			exponentially stable if, and only if 
			\begin{equation*}
			\frac{1}{2}\ln \left\vert \frac{f-1}{f+1}\right\vert =-\omega <0,
			\end{equation*}%
			which is a known result from \cite{Rideau}.
		\end{itemize}
	\end{example}
	
	\begin{remark}
		By setting 
		\begin{equation*}
		F(t)=\frac{g(\phi \circ \alpha ^{-}(t))}{g(\alpha ^{-}(t))}, \ \forall t\geq
		0,
		\end{equation*}
		with $g(t)\neq 0,$ for all $t\geq 0$, we obtain 
		\begin{equation*}
		\psi _{n}(\tau )=\prod\limits_{i=0}^{n}\left\vert F\left( \left( \alpha
		^{-}\right) ^{-1}\circ \phi ^{\left[ i\right] }(\tau )\right) \right\vert
		=\prod\limits_{i=0}^{n}\left\vert \frac{g(\phi ^{\left[ i+1\right] }(\tau )) 
		}{g(\phi ^{\left[ i\right] }(\tau )}\right\vert =\left\vert \frac{g(\phi ^{ %
				\left[ n+1\right] }(\tau ))}{g(\phi ^{\left[ 0\right] }(\tau ))}\right\vert .
		\end{equation*}
		In this case, (\ref{g}) is automatically satisfied, and since $F=\frac{1-f}{
			1+f},$ we obtain 
		\begin{equation}
		\frac{g(\alpha ^{-}(t))-g(\phi \circ \alpha ^{-}(t))}{g(\alpha
			^{-}(t))+g(\phi \circ \alpha ^{-}(t))}=f_{g}(t), \ \forall t\geq 0.
		\label{form}
		\end{equation}
		The last expression provides an explicit relation between the decay rate and
		the feedback function $f.$ This means that $f$ can be determined based on
		the desired decay rate. Formula (\ref{form}) has been used to construct $f$
		in the second and the third points in example (\ref{example}).
	\end{remark}
	
	\begin{remark}
		Examples \ref{example2} and \ref{example3} illustrate the big influence of
		the boundary curves nature on the decay rate of the solution to system (\ref%
		{stability}).
	\end{remark}
	
	\begin{remark}
		Observe that the time of extinction of the solution to system (\ref%
		{stability}) for $f\equiv 1$ is the time of exact controllability in Theorem %
		\ref{theorem control}. This can be explained by the fact that exponential
		stability implies exact controllability for time reversible systems (see for
		instance \cite[Remark 1.5]{Komornik} or \cite{Russell}). Even though our
		system is not time reversible (because of the boundary functions), we have
		seen that this implication remains true.
	\end{remark}
	
	\section{Construction of the exact solution\label{section construction}}
	
	The aim now is to find the solution $\left( p,q\right) $ to system (\ref%
	{both}) in all $Q.$ To this end, let us start by splitting $Q$ into an
	infinite number of parts. Namely 
	\begin{equation*}
	Q=\cup _{n\geq 0}\Sigma _{n}^{p}=\cup _{n\geq 0}\Sigma _{n}^{q},\text{ \ \ }%
	\Sigma _{i}^{p}\cap \Sigma _{j}^{p},\Sigma _{i}^{q}\cap \Sigma
	_{j}^{q}=\varnothing ,\text{ }i\neq j,
	\end{equation*}%
	where $\Sigma _{n}^{p},\Sigma _{n}^{q}$ are given for $n=0,1,$ by 
	\begin{eqnarray}
	\Sigma _{0}^{p} &=&\left\{ (t,x)\in Q,\text{ }t\in \lbrack 0,x)\right\} ,
	\label{region0} \\
	\Sigma _{1}^{p} &=&\left\{ (t,x)\in Q,\text{ }t-x\in \lbrack 0,\alpha
	^{-}\circ \left( \alpha ^{+}\right) ^{-1}(1))\right\} , \\
	\Sigma _{0}^{q} &=&\left\{ (t,x)\in Q,\text{ }t\in \lbrack 0,1-x)\right\} ,
	\\
	\Sigma _{1}^{q} &=&\left\{ (t,x)\in Q,\text{ }t+x\in \lbrack 1,\beta
	^{+}\circ \left( \beta ^{-}\right) ^{-1}(0))\right\} ,
	\end{eqnarray}%
	and for all $n\geq 1$ 
	\begin{equation}
	\Sigma _{2n}^{p}=\left\{ (t,x)\in Q,\text{ }t-x\in \left[ \phi ^{\left[ n-1%
		\right] }\circ \alpha ^{-}\circ \left( \alpha ^{+}\right) ^{-1}(1),\phi ^{%
		\left[ n\right] }(0)\right) \right\} ,
	\end{equation}%
	\begin{equation}
	\Sigma _{2n+1}^{p}=\left\{ (t,x)\in Q,\text{ }t-x\in \left[ \phi ^{\left[ n%
		\right] }(0),\phi ^{\left[ n\right] }\circ \alpha ^{-}\circ \left( \alpha
	^{+}\right) ^{-1}(1)\right) \right\} ,
	\end{equation}%
	\begin{equation}
	\Sigma _{2n}^{q}=\left\{ (t,x)\in Q,\text{ }t+x\in \left[ \xi ^{\left[ n-1%
		\right] }\circ \beta ^{+}\circ \left( \beta ^{-}\right) ^{-1}(0),\xi ^{\left[
		n\right] }(1)\right) \right\} ,
	\end{equation}%
	\begin{equation}
	\Sigma _{2n+1}^{q}=\left\{ (t,x)\in Q,\text{ }t+x\in \left[ \xi ^{\left[ n%
		\right] }(1),\xi ^{\left[ n\right] }\circ \beta ^{+}\circ \left( \beta
	^{-}\right) ^{-1}(0)\right) \right\} ,  \label{region1}
	\end{equation}%
	where $\xi $ is defined by 
	\begin{equation}
	\xi :=\beta ^{+}\circ \left( \beta ^{-}\right) ^{-1}\circ \alpha ^{-}\circ
	\left( \alpha ^{+}\right) ^{-1}.  \label{ksi}
	\end{equation}%
	The construction of these regions relies on the reflection of the principal
	characteristic lines with positive and negative slopes emerging from the
	points $(0,0)$ and $(0,1)$ and reflected along the boundary curves. More
	precisely, the lines $x=t$ and $x=-t+1$ emerging respectively from $(0,0)$
	and $(0,1)$ meet the curves $\left( t,\beta (t)\right) _{t\geq 0}$ and $%
	\left( t,\alpha (t)\right) _{t\geq 0}$ in the points $\left( \left( \beta
	^{-}(0)\right) ^{-1},\beta (\left( \beta ^{-}(0)\right) ^{-1})\right) $ and $%
	\left( \left( \alpha ^{+}\right) ^{-1}(1),\alpha (\left( \alpha
	^{+}(1)\right) ^{-1})\right) $ respectively. The regions $\Sigma _{0}^{p}$
	and $\Sigma _{0}^{q}$ are those located between $t=0$ and these lines. We
	can do similarly to construct the regions $\Sigma _{n}^{p},\Sigma _{n}^{q},$ 
	$n\geq 1,$ given above. In the sequel, we denote by $p_{n}$ and $q_{n}$ the
	restriction of $p$ and $q$ solutions of system (\ref{both}) on $\Sigma
	_{n}^{p}$ and $\Sigma _{n}^{q},$ $n\geq 0.$
	\begin{figure}[H]
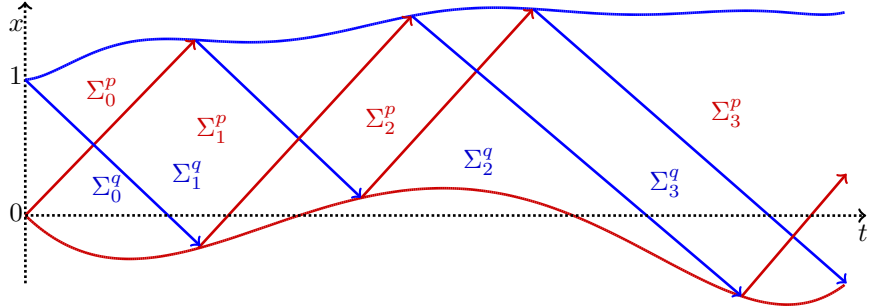

		\label{figure4} \centering
		\definecolor{qqqqcc}{rgb}{0,0,0.8} \definecolor{qqqqff}{rgb}{0,0,1}  %
		\definecolor{ccqqqq}{rgb}{0.8,0,0}  

		\caption{The regions $\Sigma^p_{i}$ are those between the red lines and $%
			\Sigma^q_{i}$ are those between the blue lines.}
	\end{figure}

	\begin{remark}
		In particular, if $\alpha \equiv 0$ and $\beta \equiv 1,$ the regions $%
		\Sigma _{n}^{p},\Sigma _{n}^{q},$ $n\geq 0,$ are simply given by 
		\begin{equation*}
		\Sigma _{n}^{p}=\left\{ (t,x)\in 
		\mathbb{R}
		_{+}\times \lbrack 0,1],\text{ }t-x\in \lbrack n-1,n)\right\} ,
		\end{equation*}%
		\begin{equation*}
		\Sigma _{n}^{q}=\left\{ (t,x)\in 
		\mathbb{R}
		_{+}\times \lbrack 0,1],\text{ }x+t\in \lbrack n,n+1)\right\} .
		\end{equation*}
	\end{remark}
	
	During the construction below, we use the standard density argument by
	assuming first that the initial states are sufficiently regular then passing
	to the limit. So, the constructed solutions must be understood in the weak
	sense. Let us start by finding $p_{0}$ and $q_{0}:$
	
	\begin{lemma}
		Let $\left( \widetilde{p},\widetilde{q}\right) \in \left[ L^{2}(0,1)\right]
		^{2}.$ The solution $(p_{0},q_{0})$ to system (\ref{both}) is given by 
		\begin{equation}
		p_{0}(t,x)=\widetilde{p}(x-t),\text{ \ }q_{0}(t,x)=\widetilde{q}(x+t).
		\label{y0+-}
		\end{equation}
	\end{lemma}
	
	\begin{proof}
		The proof readily follows from (\ref{charac}).
	\end{proof}
	
	Now, let us find the solution in the regions $\Sigma _{1}^{p},\Sigma
	_{1}^{q}:$
	
	\begin{lemma}
		Let $\left( \widetilde{p},\widetilde{q}\right) \in \left[ L^{2}(0,1)\right]
		^{2}.$ The solution $(p_{1},q_{1})$ to system (\ref{both}) is given by 
		\begin{equation}
		p_{1}(t,x)=v\left( \left( \alpha ^{-}\right) ^{-1}(t-x)\right) -F\left(
		\left( \alpha ^{-}\right) ^{-1}(t-x)\right) \widetilde{q}\left( \alpha
		^{+}\circ \left( \alpha ^{-}\right) ^{-1}(t-x)\right) ,  \label{y1-}
		\end{equation}
		\begin{equation}
		q_{1}(t,x)=-\widetilde{p}\left( -\beta ^{-}\circ \left( \beta ^{+}\right)
		^{-1}(x+t)\right) .  \label{y1+}
		\end{equation}
	\end{lemma}
	
	\begin{proof}
		By using (\ref{y0+-}), we have at the boundary curves 
		\begin{equation*}
		p_{0}(\tau ,\beta (\tau ))=\widetilde{p}\left( -\beta ^{-}(\tau )\right) , 
		\text{ \ }\tau \in \left[ 0,\left( \beta ^{-}\right) ^{-1}(0)\right) ,
		\end{equation*}
		\begin{equation*}
		q_{0}(\chi ,\alpha (\chi ))=\widetilde{q}\left( \alpha ^{+}(\chi )\right) , 
		\text{ \ \ }\chi \in \left[ 0,\left( \alpha ^{+}\right) ^{-1}(1)\right) .
		\end{equation*}
		By using the boundary conditions given in (\ref{boundary}), we get 
		\begin{eqnarray}
		&&p_{1}(\tau ,\alpha (\tau ))=v(\tau )-F(\tau )q_{0}(\tau ,\alpha (\tau ))
		\label{1} \\
		&=&v(\tau )-F(\tau )\widetilde{q}(\alpha ^{+}(\tau )),\text{ }\tau \in \left[
		0,\left( \alpha ^{+}\right) ^{-1}(1)\right) ,  \notag
		\end{eqnarray}
		\begin{equation}
		q_{1}(\chi ,\beta (\chi ))=-p_{0}(\chi ,\beta (\chi ))=-\widetilde{p}(-\beta
		^{-}(\chi )),\text{ }\chi \in \left[ 0,\left( \beta ^{-}\right)
		^{-1}(0)\right) .  \label{2}
		\end{equation}
		Consider the latter values as initial states on both regions $\Sigma
		_{1}^{p},\Sigma _{1}^{q}$ and use (\ref{charac}), we write 
		\begin{equation}
		p_{1}(t,c-t)=p_{1}(\tau ,\tau -s)\text{ \ },\text{ \ }q_{1}(\chi ,c+\chi
		)=q_{1}(\chi ,c+\chi ).  \label{3}
		\end{equation}
		By using the fact that $p$ and $q$ are constant along the characteristic
		lines $x=t-\alpha ^{-}(\tau )$ and $x=-t+\beta ^{+}(\chi )$ respectively, we
		obtain 
		\begin{equation}
		p_{1}(t,t-\alpha ^{-}(\tau ))=p_{1}(\tau ,\alpha (\tau ))=v(\tau )-F(\tau ) 
		\widetilde{q}(\alpha ^{+}(\tau )),  \label{4}
		\end{equation}
		and 
		\begin{equation}
		q_{1}(t,-t+\beta ^{+}(\chi ))=q_{1}^{-}(\chi ,\beta (\chi ))=-\widetilde{p}
		(-\beta ^{-}(\chi )).  \label{5}
		\end{equation}
		Now, letting $\left( \alpha ^{-}\right) ^{-1}(t-x)=\tau $ in (\ref{4}) and $%
		\chi =\left( \beta ^{+}\right) ^{-1}(x+t)$ in (\ref{5}) yields the desired
		result.
	\end{proof}
	
	\begin{remark}
		\label{remark}Note that $\alpha ^{+}\circ \left( \alpha ^{-}\right)
		^{-1}(t-x),$ $(t,x)\in \Sigma _{1}^{p}$ and \ $-\beta ^{-}\circ \left( \beta
		^{+}\right) ^{-1}(x+t),$ $(t,x)\in \Sigma _{1}^{q}$ belong to $(0,1)$ and
		the above expressions make perfectly sense. To clarify more things, let $%
		(t,x)\in \Sigma _{1}^{p}$ and let $\widetilde{x}(s)=s-t+x$ the line passing
		throught the point $(t,x).$ By moving backwards, this line meets the curve $%
		\left( s,\alpha (s)\right) _{s\geq 0}$ at the point $\left( \left( \alpha
		^{-}\right) ^{-1}(t-x),\alpha \left( \alpha ^{-}\right) ^{-1}(t-x)\right) $
		where $\left( \alpha ^{-}\right) ^{-1}(t-x)\in \left[ 0,\left( \alpha
		^{-}\right) ^{-1}(1)\right) .$ We use again the reflection of the
		characteristic line with negative slope passing through the latter point.
		i.e. $\widetilde{x}(s)=-s+\alpha ^{+}\circ \left( \alpha ^{-}\right)
		^{-1}(t-x)$ lying in $\Sigma _{0}^{p},$ for $s=0,$ we obtain $\widetilde{x}%
		(0)=\alpha ^{+}\circ \left( \alpha ^{-}\right) ^{-1}(t-x)$ $\in (0,1).$ We
		can do similarly for $-\beta ^{-}\circ \left( \beta ^{+}\right) ^{-1}(x+t),$ 
		$(t,x)\in \Sigma _{1}^{q}.$
	\end{remark}
	
	\begin{lemma}
		Let $\left( \widetilde{p},\widetilde{q}\right) \in \left[ L^{2}(0,1)\right]
		^{2}.$ The solution $(p_{2},q_{2})$ to system (\ref{both}) is given by 
		\begin{eqnarray}
		&&p_{2}\left( t,x\right) =v\left( \left( \alpha ^{-}\right) ^{-1}(t-x)\right)
		\label{y2-} \\
		&&+F\left( \left( \alpha ^{-}\right) ^{-1}(t-x)\right) \widetilde{p}\left(
		-\phi ^{-1}(t-x)\right) ,  \notag
		\end{eqnarray}
		\begin{eqnarray}
		&&q_{2}\left( t,x\right) =-v\left( \left( \alpha ^{-}\right) ^{-1}\circ
		\beta ^{-}\circ \left( \beta ^{+}\right) ^{-1}\left( x+t\right) \right)
		\label{y2+} \\
		&&+F\left( \left( \alpha ^{-}\right) ^{-1}\circ \beta ^{-}\circ \left( \beta
		^{+}\right) ^{-1}\left( x+t\right) \right) \widetilde{q}\left( \xi
		^{-1}\left( x+t\right) \right) ,  \notag
		\end{eqnarray}
		where $\phi $ and $\xi $ are defined in (\ref{phi}) and (\ref{ksi}).
	\end{lemma}
	\begin{proof}
		From (\ref{y1-}) and (\ref{y1+}), we have at the boundary curves 
		\begin{eqnarray}
		p_{1}(\tau ,\beta (\tau )) &=&v\left( \left( \alpha ^{-}\right) ^{-1}\circ
		\beta ^{-}(\tau )\right)   \label{6} \\
		&&-F\left( \left( \alpha ^{-}\right) ^{-1}\circ \beta ^{-}(\tau )\right) 
		\widetilde{q}\left( \alpha ^{+}\circ \left( \alpha ^{-}\right) ^{-1}\circ
		\beta ^{-}(\tau )\right) ,  \notag \\
		\tau  &\in &\left[ \left( \beta ^{-}\right) ^{-1}(0),\left( \beta
		^{-}\right) ^{-1}\circ \alpha ^{-}\circ \left( \alpha ^{+}\right)
		^{-1}(1)\right)   \notag
		\end{eqnarray}%
		and 
		\begin{eqnarray}
		q_{1}(\chi ,\alpha (\chi )) &=&-\widetilde{p}\left( -\beta ^{-}\circ \left(
		\beta ^{+}\right) ^{-1}\circ \alpha ^{+}(\chi )\right) ,  \label{7} \\
		\chi  &\in &\left[ \left( \alpha ^{+}\right) ^{-1}(1),\left( \alpha
		^{+}\right) ^{-1}\circ \beta ^{+}\circ \left( \beta ^{-}\right)
		^{-1}(0)\right) .  \notag
		\end{eqnarray}%
		In order to find $p_{2}$ and $q_{2},$ we use the boundary conditions (\ref%
		{boundary}) and the values of $p_{1}$ and $q_{1}$ at the boundary curves
		given in (\ref{6}) and (\ref{7}) as initial states. Namely, for any $\tau
		\in \left[ \left( \beta ^{-}\right) ^{-1}(0),\left( \beta ^{-}\right)
		^{-1}\circ \alpha ^{-}\circ \left( \alpha ^{+}\right) ^{-1}(1)\right) $ and 
		
		$\chi \in \left[ \left( \alpha ^{+}\right) ^{-1}(1),\left( \alpha
		^{+}\right) ^{-1}\circ \beta ^{+}\circ \left( \beta ^{-}\right)
		^{-1}(0)\right) ,$ we have along the lines \\   $x=t-\alpha ^{-}(\tau )$ and $%
		x=-t+\beta ^{+}(\chi )$ respectively 
		\begin{equation}
		p_{2}\left( t,t-\alpha ^{-}(\tau )\right) =p_{2}\left( \tau ,\alpha (\tau
		)\right) =v(\tau )-F(\tau )q_{1}\left( \tau ,\alpha (\tau )\right) ,
		\label{8}
		\end{equation}%
		\begin{equation}
		q_{2}\left( t,\beta ^{+}(\chi )-t\right) =q_{2}\left( \chi ,\beta (\chi
		)\right) =-p_{1}^{-}\left( \chi ,\beta (\chi )\right) .  \label{9}
		\end{equation}%
		Plugging (\ref{6}) and (\ref{7}) in (\ref{8}) and (\ref{9}), we get 
		\begin{equation*}
		p_{2}\left( t,t-\alpha ^{-}(\tau )\right) =v(\tau )+F(\tau )\widetilde{p}%
		\left( -\beta ^{-}\circ \left( \beta ^{+}\right) ^{-1}\circ \alpha ^{+}(\tau
		)\right) ,
		\end{equation*}%
		and 
		\begin{eqnarray*}
			&&q_{2}\left( t,\beta ^{+}(\chi )-t\right) =-v\left( \left( \alpha
			^{-}\right) ^{-1}\circ \beta ^{-}(\chi )\right)  \\
			&&+F\left( \left( \alpha ^{-}\right) ^{-1}\circ \beta ^{-}(\chi )\right) 
			\widetilde{q}\left( \alpha ^{+}\circ \left( \alpha ^{-}\right) ^{-1}\circ
			\beta ^{-}(\chi )\right) .
		\end{eqnarray*}%
		The proof follows immediately for $\tau =\left( \alpha ^{-}\right) ^{-1}(t-x)
		$ and $\left( \beta ^{+}\right) ^{-1}\left( x+t\right) =\chi $.
	\end{proof}
	
	\begin{remark}
		In the same spirit of Remark \ref{remark}, the expressions (\ref{y2-}) and (%
		\ref{y2+}) make perfectly sense. We can use the same reasonning to show that 
		\begin{eqnarray*}
			-\beta ^{-}\circ \left( \beta ^{+}\right) ^{-1}\circ \alpha ^{+}\circ \left(
			\alpha ^{-}\right) ^{-1}(t-x) &\in &(0,1),\text{ }\forall (t,x)\in \Sigma
			_{2}^{p}, \\
			\alpha ^{+}\circ \left( \alpha ^{-}\right) ^{-1}\circ \beta ^{-}\circ \left(
			\beta ^{+}\right) ^{-1}\left( x+t\right) &\in &(0,1),\text{ }\forall
			(t,x)\in \Sigma _{2}^{q}.
		\end{eqnarray*}%
		More generally, we have:
	\end{remark}
	
	\begin{lemma}
		Let $\left( \widetilde{p},\widetilde{q}\right) \in \left[ L^{2}(0,1)\right]
		^{2}.$ The solutions $p_{2n+1},p_{2n+2},q_{2n+1},q_{2n+2},$ $n\geq 1,$ to
		system (\ref{both}) are given by 
		\begin{eqnarray}
		&&p_{2n+1}(t,x)  \label{p2n+1} \\
		&=&\sum_{k=0}^{n}v\left( \left( \alpha ^{-}\right) ^{-1}\circ \left( \phi
		^{-1}\right) ^{[k]}(t-x)\right) \prod\limits_{i=0}^{k-1}F\left( \left(
		\alpha ^{-}\right) ^{-1}\circ \left( \phi ^{-1}\right) ^{[i]}(t-x)\right) 
		\notag \\
		&&-\widetilde{q}\left( \left( \xi ^{-1}\right) ^{\left[ n\right] }\circ
		\alpha ^{+}\circ \left( \alpha ^{-}\right) ^{-1}(t-x)\right)
		\prod\limits_{k=0}^{n}F\left( \left( \alpha ^{-}\right) ^{-1}\circ \left(
		\phi ^{-1}\right) ^{[k]}(t-x)\right) ,  \notag
		\end{eqnarray}%
		\begin{eqnarray}
		&&p_{2n+2}\left( t,x\right)  \label{p2n+2} \\
		&=&\sum_{k=0}^{n}v\left( \left( \alpha ^{-}\right) ^{-1}\circ \left( \phi
		^{-1}\right) ^{[k]}(t-x)\right) \prod\limits_{i=0}^{k-1}F\left( \left(
		\alpha ^{-}\right) ^{-1}\circ \left( \phi ^{-1}\right) ^{[i]}(t-x)\right) 
		\notag \\
		&&+\widetilde{p}\left( -\left( \phi ^{-1}\right) ^{\left[ n+1\right] }\circ
		(t-x)\right) \prod\limits_{k=0}^{n}F\left( \left( \alpha ^{-}\right)
		^{-1}\circ \left( \phi ^{-1}\right) ^{[k]}(t-x)\right) ,  \notag
		\end{eqnarray}%
		\begin{eqnarray}
		&&q_{2n+1}(t,x)  \label{q2n+1} \\
		&=&-\sum_{k=0}^{n-1}v\left( \left( \alpha ^{-}\right) ^{-1}\circ \left( \phi
		^{-1}\right) ^{[k]}\circ \beta ^{-}\circ \left( \beta ^{+}\right)
		^{-1}\left( x+t\right) \right) \times  \notag \\
		&&\prod\limits_{i=0}^{k-1}F\left( \left( \alpha ^{-}\right) ^{-1}\circ
		\left( \phi ^{-1}\right) ^{[i]}\circ \beta ^{-}\circ \left( \beta
		^{+}\right) ^{-1}\left( x+t\right) \right)  \notag \\
		&&-\widetilde{p}\left( -\left( \phi ^{-1}\right) ^{\left[ n\right] }\circ
		\beta ^{-}\circ \left( \beta ^{+}\right) ^{-1}\left( x+t\right) \right)
		\times  \notag \\
		&&\prod\limits_{k=0}^{n-1}F\left( \left( \alpha ^{-}\right) ^{-1}\circ
		\left( \phi ^{-1}\right) ^{[k]}\circ \beta ^{-}\circ \left( \beta
		^{+}\right) ^{-1}\left( x+t\right) \right) ,  \notag
		\end{eqnarray}%
		\begin{eqnarray}
		&&q_{2n+2}(t,x)  \label{q2n+2} \\
		&=&-\sum_{k=0}^{n}v\left( \left( \alpha ^{-}\right) ^{-1}\circ \left( \phi
		^{-1}\right) ^{[k]}\circ \beta ^{-}\circ \left( \beta ^{+}\right)
		^{-1}(x+t)\right) \times  \notag \\
		&&\prod\limits_{i=0}^{k-1}F\left( \left( \alpha ^{-}\right) ^{-1}\circ
		\left( \phi ^{-1}\right) ^{[i]}\circ \beta ^{-}\circ \left( \beta
		^{+}\right) ^{-1}(x+t)\right)  \notag \\
		&&+\widetilde{q}\left( \left( \xi ^{-1}\right) ^{\left[ n+1\right]
		}(x+t)\right) \times  \notag \\
		&&\prod\limits_{k=0}^{n}F\left( \left( \alpha ^{-}\right) ^{-1}\circ \left(
		\phi ^{-1}\right) ^{[k]}\circ \beta ^{-}\circ \left( \beta ^{+}\right)
		^{-1}(x+t)\right) ,  \notag
		\end{eqnarray}%
		with the convention $\prod\limits_{k=0}^{-1}=1.$ The functions $\phi $ and $%
		\xi $ are defined in (\ref{phi}) and (\ref{ksi}).
	\end{lemma}
	
	\begin{proof}
		The above expressions can be proved by induction. Let us start by proving (%
		\ref{q2n+2}). At the boundary $x=\beta (t)$, (\ref{p2n+1}) becomes 
		\begin{eqnarray*}
			&&p_{2n+1}(t,\beta (t)) \\
			&=&\sum_{k=0}^{n}v\left( \left( \alpha ^{-}\right) ^{-1}\circ \left( \phi
			^{-1}\right) ^{[k]}\circ \beta ^{-}(t)\right)
			\prod\limits_{i=0}^{k-1}F\left( \left( \alpha ^{-}\right) ^{-1}\circ \left(
			\phi ^{-1}\right) ^{[i]}\circ \beta ^{-}(t)\right) \\
			&&-\widetilde{q}\left( \left( \xi ^{-1}\right) ^{\left[ n\right] }\circ
			\alpha ^{+}\circ \left( \alpha ^{-}\right) ^{-1}\circ \beta ^{-}(t)\right)
			\prod\limits_{k=0}^{n}F\left( \left( \alpha ^{-}\right) ^{-1}\circ \left(
			\phi ^{-1}\right) ^{[k]}\circ \beta ^{-}(t)\right) .
		\end{eqnarray*}%
		Now, we use the boundary condition given in (\ref{boundary}), i.e. 
		\begin{eqnarray*}
			q_{2n+2}(\chi ,\beta (\chi )) &=&-p_{2n+1}(\chi ,\beta (\chi )),\text{ \ } \\
			\chi &{\small \in }&\left[ \text{ }\left( \beta ^{-1}\right) ^{-1}\circ 
			\text{ }\phi ^{\left[ n\right] }(0),\text{ }\left( \beta ^{-1}\right)
			^{-1}\circ \phi ^{\left[ n\right] }\circ \alpha ^{-}\circ \left( \alpha
			^{+}\right) ^{-1}(1)\right) ,
		\end{eqnarray*}%
		we find 
		\begin{eqnarray}
		&&q_{2n+2}(\chi ,\beta (\chi ))  \label{s} \\
		&=&-\sum_{k=0}^{n}v\left( \left( \alpha ^{-}\right) ^{-1}\circ \left( \phi
		^{-1}\right) ^{[k]}\circ \beta ^{-}(\chi )\right)
		\prod\limits_{i=0}^{k-1}F\left( \left( \alpha ^{-}\right) ^{-1}\circ \left(
		\phi ^{-1}\right) ^{[i]}\circ \beta ^{-}(t)\right)  \notag \\
		&&+\widetilde{q}\left( \left( \xi ^{-1}\right) ^{\left[ n\right] }\circ
		\alpha ^{+}\circ \left( \alpha ^{-}\right) ^{-1}\circ \beta ^{-}(\chi
		)\right) \prod\limits_{k=0}^{n}F\left( \left( \alpha ^{-}\right) ^{-1}\circ
		\left( \phi ^{-1}\right) ^{[k]}\circ \beta ^{-}(\chi )\right) .  \notag
		\end{eqnarray}%
		Since $q$ is constant along the characteristic lines of the form $x=c-t$, in
		particular, on the line $x=\beta ^{+}(\chi )-t,$ we have 
		\begin{equation*}
		q_{2n+2}(t,\beta ^{+}(\chi )-t)=q_{2n+2}(\chi ,\beta (\chi )).
		\end{equation*}%
		Finally, by letting $\chi =\left( \beta ^{+}\right) ^{-1}(x+t)$ in (\ref{s}%
		), we obtain the formula in (\ref{q2n+2}). Let us do similarly for $%
		p_{2n+2}. $ By taking (\ref{q2n+1}) for $x=\alpha (t),$ we obtain 
		\begin{eqnarray}
		&&q_{2n+1}(t,\alpha (t))  \label{z} \\
		&=&-\sum_{k=0}^{n-1}v\left( \left( \alpha ^{-}\right) ^{-1}\circ \left( \phi
		^{-1}\right) ^{[k]}\circ \beta ^{-}\circ \left( \beta ^{+}\right) ^{-1}\circ
		\alpha ^{+}(t)\right) \times  \notag \\
		&&\prod\limits_{i=0}^{k-1}F\left( \left( \alpha ^{-}\right) ^{-1}\circ
		\left( \phi ^{-1}\right) ^{[i]}\circ \beta ^{-}\circ \left( \beta
		^{+}\right) ^{-1}\circ \alpha ^{+}(t)\right)  \notag \\
		&&-\widetilde{p}\left( -\left( \phi ^{-1}\right) ^{\left[ n\right] }\circ
		\beta ^{-}\circ \left( \beta ^{+}\right) ^{-1}\circ \alpha ^{+}(t)\right)
		\times  \notag \\
		&&\prod\limits_{k=0}^{n-1}F\left( \left( \alpha ^{-}\right) ^{-1}\circ
		\left( \phi ^{-1}\right) ^{[k]}\circ \beta ^{-}\circ \left( \beta
		^{+}\right) ^{-1}\circ \alpha ^{+}(t)\right) .  \notag
		\end{eqnarray}%
		Using the boundary condition 
		\begin{eqnarray*}
			p_{2n+2}(\tau ,\alpha (\tau )) &=&v(\tau )-F(\tau )q_{2n+1}(\tau ,\alpha
			(\tau )),\text{ \ } \\
			\tau &\in &\left[ \left( \alpha ^{+}\right) ^{-1}\circ \xi ^{\left[ n\right]
			}(1),\left( \alpha ^{+}\right) ^{-1}\circ \xi ^{\left[ n\right] }\circ \beta
			^{+}\circ \left( \beta ^{-}\right) ^{-1}(0)\right) ,
		\end{eqnarray*}%
		and the fact that $q$ is constant along the characteristic lines $x=c-t$, in
		particular, on the line $x=t-\alpha ^{-}(\tau ),$ we obtain 
		\begin{equation}
		p_{2n+2}(\tau ,t-\alpha ^{-}(\tau ))=v(\tau )-F(\tau )q_{2n+1}(\tau ,\alpha
		(\tau )).  \label{h}
		\end{equation}%
		By letting $\tau =\left( \alpha ^{-}\right) ^{-1}(t-x)$ in (\ref{z}) and
		plugging the result in (\ref{h}) then using the definition of $\phi $ given
		in (\ref{phi}), we get 
		\begin{eqnarray*}
			&&p_{2n+2}(t,x) \\
			&=&v(\left( \alpha ^{-}\right) ^{-1}(t-x))+\sum_{k=0}^{n-1}v\left( \left(
			\alpha ^{-}\right) ^{-1}\circ \left( \phi ^{-1}\right) ^{[k+1]}(t-x)\right)
			\times \\
			&&\prod\limits_{i=0}^{k-1}F(\left( \alpha ^{-}\right) ^{-1}(t-x))F\left(
			\left( \alpha ^{-}\right) ^{-1}\circ \left( \phi ^{-1}\right)
			^{[i+1]}(t-x)\right) \\
			&&+F(\left( \alpha ^{-}\right) ^{-1}(t-x))\widetilde{p}\left( -\left( \phi
			^{-1}\right) ^{\left[ n+1\right] }(t-x)\right) \times \\
			&&\prod\limits_{k=0}^{n-1}F\left( \left( \alpha ^{-}\right) ^{-1}\circ
			\left( \phi ^{-1}\right) ^{[k+1]}(t-x)\right) .
		\end{eqnarray*}%
		After some manipulation we obtain the formula in (\ref{p2n+2}).
	\end{proof}
	
	\begin{remark}
		From what preceed, it is not difficult to see that the solution $(p,q)$ to
		system (\ref{both}) satisfies the regularity given in (\ref{regularity}).
	\end{remark}
	
	\begin{remark}
		More generally, if $\left( \widetilde{p},\widetilde{q},v,F\right) \in \left[
		L^{\theta }(0,1)\right] ^{2}\times L_{\mathrm{loc}}^{\theta }(0,\infty
		)\times L^{\eta }(0,\infty ),$ we can see from (\ref{p2n+1})-(\ref{q2n+2})
		that the solution $(p,q)$ to system (\ref{both}) satisfies the regularity 
		\begin{equation*}
		(p,q)\in \mathcal{C}(0,t;\left[ L^{r}(\alpha (t),\beta (t))\right] ^{2}),%
		\text{ }t\geq 0,
		\end{equation*}%
		with $\frac{1}{\theta }+\frac{1}{\eta }=\frac{1}{r},$ $\theta ,\eta \in
		\lbrack 1,\infty ).$
	\end{remark}
	
	\section{Proof of main results}
	
	\subsection{Proof of the controllability theorem}
	
	Let $F\equiv 1$ in (\ref{y1-}),(\ref{y1+}),(\ref{y2-}) and (\ref{y2+}). The
	solution $p_{1}$ sees the control immediately for $t\geq 0$, on the
	contrary, the component $q_{1}$ has to wait one more reflection on the curve 
	$(t,\alpha (t))_{t\geq 0}$ to see it as soon as $\ t\geq \left( \beta
	^{-}\right) ^{-1}(0)$. Le us start by proving the necessary part:
	
	\begin{proposition}
		If $T<$ $T^{\ast }=\left( \alpha ^{+}\right) ^{-1}\circ \beta ^{+}\circ
		\left( \beta ^{-}\right) ^{-1}(0),$ then system (\ref{control2}) is not
		exactly controllable at time $T$.
	\end{proposition}
	
	\begin{proof}
		To prove this lemma, we make use of the expressions of the exact solution
		given in (\ref{y1+}) and (\ref{y2+}). Let $T_{\varepsilon }^{\ast }=T^{\ast
		}-\varepsilon $ for sufficiently small $\varepsilon >0$; the solution $q$ at
		this time is given by 
		\begin{equation*}
		q(T_{\varepsilon }^{\ast },x)=\left\{ 
		\begin{array}{ccc}
		q_{1}^{+}(T_{\varepsilon }^{\ast },x) & \mathrm{if} & x\in \left[ \alpha
		(T_{\varepsilon }^{\ast }),T_{\varepsilon }^{\ast }-\beta ^{+}\circ \left(
		\beta ^{-}\right) ^{-1}(0)\right) , \\ 
		q_{2}^{+}(T_{\varepsilon }^{\ast },x) & \mathrm{if} & x\in \left[
		T_{\varepsilon }^{\ast }-\beta ^{+}\circ \left( \beta ^{-}\right)
		^{-1}(0),\beta (T_{\varepsilon }^{\ast })\right) .%
		\end{array}
		\right.
		\end{equation*}
		Thus, system (\ref{control}) will be never exactly controllable since we
		have for any initial state $\widetilde{p}$ and any target state $k$ 
		\begin{equation*}
		q(T_{\varepsilon }^{\ast },x)=-\widetilde{p}\left( -\beta ^{-}\circ \left(
		\beta ^{+}\right) ^{-1}(x+T_{\varepsilon })\right) =k(x),\text{ }x\in \left[
		\alpha (T_{\varepsilon }^{\ast }),T_{\varepsilon }^{\ast }-\beta ^{+}\circ
		\left( \beta ^{-}\right) ^{-1}(0)\right) ,
		\end{equation*}
		which is clearly a violating of the initial states.
	\end{proof}
	
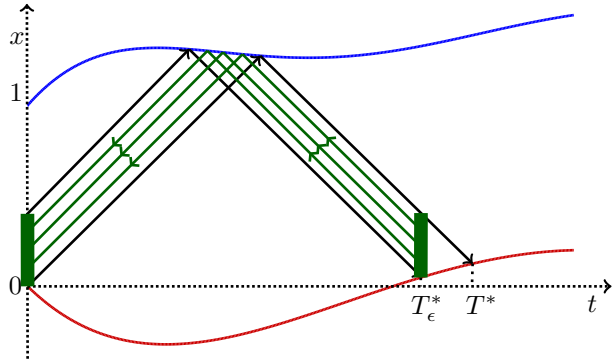
\begin{figure}[H]
		\centering
		\definecolor{qqwuqq}{rgb}{0,0.39215686274509803,0}  %
		\definecolor{ccqqqq}{rgb}{0.8,0,0} \definecolor{qqqqff}{rgb}{0,0,1}  

	\end{figure}
	
	Now, we prove the sufficient part:
	
	\begin{proposition}
		If $T\geq $ $T^{\ast }=\left( \alpha ^{+}\right) ^{-1}\circ \beta ^{+}\circ
		\left( \beta ^{-}\right) ^{-1}(0),$ then system (\ref{both}) is exactly
		controllable at time $T$.
	\end{proposition}
	
	\begin{proof}
		It suffices to prove it for $T=$ $T^{\ast }$. Let $\left( h,k\right) \in
		L^{2}(\alpha (T^{\ast }),\beta (T^{\ast }))$ be a target state and let $%
		T^{\ast \ast }=\left( \beta ^{-}\right) ^{-1}\circ \alpha ^{-}\circ \left(
		\alpha ^{+}\right) ^{-1}(1).$ We have three possible configurations:
		
		\textbf{Case 1}: $T^{\ast \ast }=T^{\ast }$
		
		In this case, we have $p(T^{\ast })=p_{2}(T^{\ast })$ and $q(T^{\ast
		})=q_{2}(T^{\ast }),$ then by making use of (\ref{y2-}) and (\ref{y2+}) we
		obtain 
		\begin{equation*}
		h\left( x\right) =p_{2}\left( T^{\ast },x\right) =v\left( \left( \alpha
		^{-}\right) ^{-1}(T^{\ast }-x)\right) +\widetilde{p}\left( -\phi
		^{-1}(T^{\ast }-x)\right) ,\text{ \ }x\in (\alpha (T^{\ast }),\beta (T^{\ast
		}))
		\end{equation*}%
		\begin{equation*}
		k\left( x\right) =-v\left( \left( \alpha ^{-}\right) ^{-1}\circ \beta
		^{-}\circ \left( \beta ^{+}\right) ^{-1}\left( x+T^{\ast }\right) \right) +%
		\widetilde{q}\left( \xi ^{-1}\left( x+T^{\ast }\right) \right) ,\text{ }x\in
		(\alpha (T^{\ast }),\beta (T^{\ast })).
		\end{equation*}%
		Therefore, the control $v$ is given by 
		\begin{equation*}
		v(t)=\left\{ 
		\begin{array}{ccc}
		h\left( T^{\ast }-\alpha ^{-}(t)\right) -\widetilde{p}\left( -\phi
		^{-1}\circ \alpha ^{-}(t)\right) , & \mathrm{if} & t\in \left( \left( \alpha
		^{-}\right) ^{-1}\circ \beta ^{-}(T^{\ast }),T^{\ast }\right) , \\ 
		&  &  \\ 
		\begin{array}{c}
		\widetilde{q}\left( \alpha ^{+}(t)\right) \\ 
		-k\left( \beta ^{+}\circ \left( \beta ^{-}\right) ^{-1}\circ \alpha
		^{-}(t)-T^{\ast }\right) ,%
		\end{array}
		& \mathrm{if} & \text{ }t\in \left( 0,\left( \alpha ^{-}\right) ^{-1}\circ
		\beta ^{-}(T^{\ast })\right) .\text{ \ \ }%
		\end{array}%
		\right.
		\end{equation*}
		
		\textbf{Case 2}: $T^{\ast \ast }<T^{\ast }$
		
		In this case, $p(T^{\ast })$ and $q(T^{\ast })$ are defined by 
		\begin{equation*}
		p(T^{\ast },x)=\left\{ 
		\begin{array}{ccc}
		p_{1}(T^{\ast },x), & \mathrm{if} & x\in \left( T^{\ast }-\alpha ^{-}\circ
		\left( \alpha ^{+}\right) ^{-1}(1),\beta (T^{\ast })\right) , \\ 
		p_{2}(T^{\ast },x), & \mathrm{if} & x\in \left( \alpha (T^{\ast }),T^{\ast
		}-\alpha ^{-}\circ \left( \alpha ^{+}\right) ^{-1}(1)\right) ,%
		\end{array}%
		\right.
		\end{equation*}%
		and $q(T^{\ast })=q_{2}(T^{\ast }).$ Thus, by making use of (\ref{y1-}),(\ref%
		{y2-}) and (\ref{y2+}), then making some variable substitutions, we arrive
		at 
		\begin{equation*}
		v(t)=\left\{ 
		\begin{array}{ccc}
		h_{1}(T^{\ast }-\alpha ^{-}(t))+\widetilde{q}\left( \alpha ^{+}(t)\right) ,
		& \mathrm{if} & t\in \left( \left( \alpha ^{-}\right) ^{-1}\circ \beta
		^{-}(T^{\ast }),\left( \alpha ^{+}\right) ^{-1}(1)\right) , \\ 
		&  &  \\ 
		h_{2}\left( T^{\ast }-\alpha ^{-}(t)\right) -\widetilde{p}\left( -\phi
		^{-1}\circ \alpha ^{-}(t)\right) , & \mathrm{if} & t\in \left( \left( \alpha
		^{+}\right) ^{-1}(1),T^{\ast }\right) ,\text{ \ \ \ \ \ \ \ \ \ \ \ \ \ \ \
			\ \ } \\ 
		&  &  \\ 
		\begin{array}{c}
		\widetilde{q}\left( \alpha ^{+}(t)\right) \\ 
		-k\left( \beta ^{+}\circ \left( \beta ^{-}\right) ^{-1}\circ \alpha
		^{-}(t)-T^{\ast }\right) ,%
		\end{array}
		& \mathrm{if} & t\in \left( 0,\left( \alpha ^{-}\right) ^{-1}\circ \beta
		^{-}(T^{\ast })\right) ,\text{ \ \ \ \ \ \ \ \ \ \ \ }%
		\end{array}%
		\right.
		\end{equation*}%
		where $h_{1}$ and $h_{2}$ are the restrictions of the target state $h$ on
		the regions $\Sigma _{1}^{p}$ and $\Sigma _{2}^{p}$ respectively.
		
		\textbf{Case 3}: $T^{\ast \ast }>T^{\ast}$
		
		In this case, we have $p(T^{\ast })=p_{2}(T^{\ast })$, and $q(T^{\ast })$ is
		defined by 
		\begin{equation*}
		q(T^{\ast },x)=\left\{ 
		\begin{array}{ccc}
		q_{2}(T^{\ast },x), & \mathrm{if} & x\in \left( \alpha \left( T^{\ast
		}\right) ,\xi (1)-T^{\ast }\right) , \\ 
		q_{3}(T^{\ast },x), & \mathrm{if} & x\in \left( \xi (1)-T^{\ast },\beta
		(T^{\ast })\right) \ .%
		\end{array}%
		\right.
		\end{equation*}%
		By using (\ref{y2-}), (\ref{y2+}) and (\ref{q2n+1}) for $n=1$ and $t=T^{\ast
		},$ then making some variable substitutions, we obtain 
		\begin{equation*}
		v(t)=\left\{ 
		\begin{array}{ccc}
		h\left( T^{\ast }-\alpha ^{-}(t)\right) -\widetilde{p}\left( -\phi
		^{-1}\circ \alpha ^{-}(t)\right) , & \mathrm{if} & t\in \left( \left( \alpha
		^{-}\right) ^{-1}\circ \beta ^{-}(T^{\ast }),T^{\ast }\right) ,\text{ \ \ \
			\ \ \ \ \ \ } \\ 
		&  &  \\ 
		\begin{array}{c}
		\widetilde{q}\left( \alpha ^{+}(t)\right) \\ 
		-k_{2}(\beta ^{+}\circ \left( \beta ^{-}\right) ^{-1}\circ \alpha
		^{-}(t)-T^{\ast }),%
		\end{array}
		& \mathrm{if} & t\in \left( 0,\left( \alpha ^{+}\right) ^{-1}(1)\right) ,%
		\text{ \ \ \ \ \ \ \ \ \ \ \ \ \ \ \ \ \ \ \ } \\ 
		&  &  \\ 
		\begin{array}{c}
		-\widetilde{p}\left( -\beta ^{-}\circ \left( \beta ^{+}\right) ^{-1}\circ
		\alpha ^{+}(t)\right) \\ 
		-k_{3}(\beta ^{+}\circ \left( \beta ^{-}\right) ^{-1}\circ \alpha ^{-}\left(
		t\right) -T^{\ast }),%
		\end{array}
		& \mathrm{if} & t\in \left( \left( \alpha ^{+}\right) ^{-1}(1),\left( \alpha
		^{-}\right) ^{-1}\circ \beta ^{-}(T^{\ast })\right) ,%
		\end{array}%
		\right.
		\end{equation*}%
		where $k_{2}$ and $k_{3}$ are the restrictions of the target state $k$ on
		the regions $\Sigma _{2}^{q}$ and $\Sigma _{3}^{q}$ respectively. The above
		expressions are well defined and the control $v$ is uniquely determined on $%
		[0,T^{\ast })$. In particular, from (\ref{y2-}) and (\ref{y2+}), we can see
		that the control 
		\begin{equation*}
		v(t)=\left\{ 
		\begin{array}{ccc}
		\widetilde{q}\left( \alpha ^{+}(t)\right) , & \mathrm{if} & t\in \left[
		0,\left( \alpha ^{+}\right) ^{-1}(1)\right) ,\text{{}} \\ 
		-\widetilde{p}\left( -\beta ^{-}\circ \left( \beta ^{+}\right) ^{-1}\circ
		\alpha ^{+}(t)\right) , & \mathrm{if} & \ t\in \left[ \left( \alpha
		^{+}\right) ^{-1}(1),T^{\ast }\right) , \\ 
		0, & \mathrm{if} & t\geq T^{\ast },\text{ \ \ \ \ \ \ \ \ \ \ \ \ \ \ \ \ }%
		\end{array}%
		\right.
		\end{equation*}%
		makes $p_{2}$ and $q_{2}$ vanish, then by the boundary conditions given in (%
		\ref{boundary}) all the solutions $p_{n},q_{n},$ $n\geq 2$, will be zero. To
		get an explicit formula of the control $u,$ it suffices to inverse the
		transformation defined in (\ref{Riemann}), then using the compatibility
		condition $y_{0}(0)=u(0)$ to obtain (\ref{co}).
	\end{proof}
	
	\begin{remark}
		Since we have an explicit formula of the solution for all $t\geq 0,$ we can
		prove that exact controllability holds at any time $T>T^{\ast }$ with loss
		of uniqueness of the control.
	\end{remark}
	
	\subsection{Proof of the stability theorem}
	
	In this subsection, we let $v\equiv 0$. We start by proving the sufficient
	part.
	
	At time $t\geq 0,$ the components $p(t)$ and $q(t)$ might involve at most
	three values of the restrictive solutions $p_{n}(t)$ and $q_{n}(t)$
	respectively on the contrary of the cylindrical case where $p(t)$ and $q(t)$
	might involve at most two values (see Figure \ref{figure4}), (if $p(t)$ or $%
	q(t)$ are defined on four regions, we obtain $\alpha (t)>\beta (t)$). More
	precisely, we have for the component $p$:\newline
	\textbf{Case 1}: $t\in \left[ \left( \alpha ^{-}\right) ^{-1}\circ \phi ^{%
		\left[ n-1\right] }\circ \alpha ^{-}\circ \left( \alpha ^{+}\right)
	^{-1}(1),\left( \alpha ^{-}\right) ^{-1}\circ \phi ^{\left[ n\right]
	}(0)\right) .$
	
	In this case, $p(t)$ might expressed in function of $%
	p_{2n-1}(t),p_{2n}(t),p_{2n+1}(t),$ 
	\begin{equation}
	p(t,x)=\left\{ 
	\begin{array}{ccc}
	p_{2n-1}(t,x), & \text{\textrm{if}} & x\in I_{1}(t):=\left[ t-\phi ^{\left[
		n-1\right] }\circ \alpha ^{-}\circ \left( \alpha ^{+}\right) ^{-1}(1),\beta
	(t)\right) ,\text{ \ \ \ \ \ \ \ } \\ 
	p_{2n}(t,x),\text{ \ } & \text{\textrm{if}} & x\in I_{2}(t):=\left[ t-\phi ^{%
		\left[ n\right] }(0),t-\phi ^{\left[ n-1\right] }\circ \alpha ^{-}\circ
	\left( \alpha ^{+}\right) ^{-1}(1)\right) , \\ 
	p_{2n+1}(t,x), & \text{\textrm{if}} & x\in I_{3}(t):=\left[ \alpha
	(t),t-\phi ^{\left[ n\right] }(0)\right) .\text{ \ \ \ \ \ \ \ \ \ \ \ \ \ \
		\ \ \ \ \ \ \ \ \ \ \ \ \ \ \ \ }%
	\end{array}%
	\right.  \label{p1}
	\end{equation}%
	By definition of the regions $\Sigma _{n}^{p},$ $n\geq 0,$ given in (\ref%
	{region0})-(\ref{region1}), we have for $k=1,2,3$ 
	\begin{eqnarray*}
		&&\left\{ (t,x)\in \left[ \left( \alpha ^{-}\right) ^{-1}\circ \phi ^{\left[
			n-1\right] }\circ \alpha ^{-}\circ \left( \alpha ^{+}\right) ^{-1}(1),\left(
		\alpha ^{-}\right) ^{-1}\circ \phi ^{\left[ n\right] }(0)\right) \times
		I_{k}(t)\right\} \\
		&\subset &\Sigma _{2n+k-2}^{p}.
	\end{eqnarray*}%
	Consequently, 
	\begin{eqnarray}
	\left\Vert p(t)\right\Vert _{L^{2}(\alpha (t),\beta (t))}^{2}
	&=&\sum_{k=1}^{3}\left\Vert p_{2n+k-2}(t)\right\Vert _{L^{2}(I_{k}(t))}^{2}
	\label{side} \\
	&\leq &\sum_{k=1}^{3}\int_{(t,x)\in \Sigma _{2n+k-2}^{p}}\left\vert
	p_{2n+k-2}(t,x)\right\vert ^{2}dx,  \notag
	\end{eqnarray}%
	which leads us to estimate the right hand side of (\ref{side}). By using the
	exact solution formulas given in (\ref{p2n+1}) and (\ref{p2n+2}), we obtain
	for $k=1,2,3$ 
	\begin{eqnarray*}
		&&\sum_{k=1}^{3}\int_{(t,x)\in \Sigma _{2n+k-2}^{p}}\left\vert
		p_{2n+k-2}(t,x)\right\vert ^{2}dx \\
		&\leq &\left\Vert \left( \widetilde{p},\widetilde{q}\right) \right\Vert
		_{L^{2}(0,1)}^{2}\sum_{k=1}^{3}\underset{x,(t,x)\in \Sigma _{2n+k-2}^{p}}{%
			\sup }\prod\limits_{i=0}^{n-1+\left[ \frac{k-1}{2}\right] }\left\vert
		F\left( \left( \alpha ^{-}\right) ^{-1}\circ \left( \phi ^{-1}\right)
		^{[i]}(t-x)\right) \right\vert .
	\end{eqnarray*}%
	By definition of the regions $\Sigma _{n}^{p},$ $n\geq 0$ given in (\ref%
	{region0})-(\ref{region1}), we have 
	\begin{eqnarray*}
		(t,x) &\in &\Sigma _{2n}^{p}\text{ \ \ }\Leftrightarrow t-x\in \left[ \phi ^{%
			\left[ n-1\right] }\circ \alpha ^{-}\circ \left( \alpha ^{+}\right)
		^{-1}(1),\phi ^{\left[ n\right] }(0\right) , \\
		(t,x) &\in &\Sigma _{2n+1}^{p}\Leftrightarrow t-x\in \left[ \phi ^{\left[ n%
			\right] }(0),\phi ^{\left[ n\right] }\circ \alpha ^{-}\circ \left( \alpha
		^{+}\right) ^{-1}(1)\right) ,
	\end{eqnarray*}%
	therefore, there exist a sequences $\tau _{1}^{n}(t,x)\in \left[ \alpha
	^{-}\circ \left( \alpha ^{+}\right) ^{-1}(1),\phi (0)\right) $ and $\tau
	_{2}^{n}(t,x)\in \left[ 0,\alpha ^{-}\circ \left( \alpha ^{+}\right)
	^{-1}(1)\right) $ such that 
	\begin{eqnarray}
	(t,x) &\in &\Sigma _{2n}^{p}\text{ \ \ }\Leftrightarrow t-x=\phi ^{\left[ n-1%
		\right] }\left( \tau _{1}^{n}(t,x)\right) ,  \label{pppp} \\
	(t,x) &\in &\Sigma _{2n+1}^{p}\Leftrightarrow t-x=\phi ^{\left[ n\right]
	}(\tau _{2}^{n}(t,x)).  \notag
	\end{eqnarray}%
	Observe that when $(t,x)$ runs $\Sigma _{2n}^{p}$ (resp. $\Sigma _{2n+1}^{p}$%
	), the bounded sequence $\tau _{1}^{n}(t,x)$ (resp. $\tau _{2}^{n}(t,x)$)
	rises $\left[ \alpha ^{-}\circ \left( \alpha ^{+}\right) ^{-1}(1),\phi
	(0)\right) $ (resp. $\left[ 0,\alpha ^{-}\circ \left( \alpha ^{+}\right)
	^{-1}(1)\right) $). These sequences will play the role of two parameters $%
	\tau _{1}\in \left[ \alpha ^{-}\circ \left( \alpha ^{+}\right) ^{-1}(1),\phi
	(0)\right) $ and $\tau _{2}\in \left[ 0,\alpha ^{-}\circ \left( \alpha
	^{+}\right) ^{-1}(1)\right) $). With these notations, we have 
	\begin{eqnarray*}
		&&\sum_{k=1}^{3}\underset{x,(t,x)\in \Sigma _{2n+k-2}^{p}}{\sup }%
		\prod\limits_{i=0}^{n+k-2}\left\vert F\left( \left( \alpha ^{-}\right)
		^{-1}\circ \left( \phi ^{-1}\right) ^{[i]}(t-x)\right) \right\vert \\
		&\leq &\underset{\tau _{2}\in \left[ 0,\alpha ^{-}\circ \left( \alpha
			^{+}\right) ^{-1}(1)\right) }{\sup }\prod\limits_{i=0}^{n-1}\left\vert
		F\left( \left( \alpha ^{-}\right) ^{-1}\circ \left( \phi ^{-1}\right)
		^{[i]}\circ \phi ^{\left[ n-1\right] }(\tau _{1})\right) \right\vert \\
		&&+\underset{\tau _{1}\in \left[ \alpha ^{-}\circ \left( \alpha ^{+}\right)
			^{-1}(1),\phi (0)\right) }{\sup }\prod\limits_{i=0}^{n-1}\left\vert F\left(
		\left( \alpha ^{-}\right) ^{-1}\circ \left( \phi ^{-1}\right) ^{[i]}\circ
		\phi ^{\left[ n-1\right] }(\tau _{2})\right) \right\vert \\
		&&+\underset{\tau _{2}\in \left[ 0,\alpha ^{-}\circ \left( \alpha
			^{+}\right) ^{-1}(1)\right) }{\sup }\prod\limits_{i=0}^{n}\left\vert F\left(
		\left( \alpha ^{-}\right) ^{-1}\circ \left( \phi ^{-1}\right) ^{[i]}\circ
		\phi ^{\left[ n\right] }(\tau _{1})\right) \right\vert \\
		&=&\underset{\tau \in \left[ 0,\alpha ^{-}\circ \left( \alpha ^{+}\right)
			^{-1}(1)\right) }{\sup }\psi _{n-1}(s)+\underset{\tau \in \left[ \alpha
			^{-}\circ \left( \alpha ^{+}\right) ^{-1}(1),\phi (0)\right) }{\sup }\psi
		_{n-1}(\tau ) \\
		&&+\underset{\tau \in \left[ 0,\alpha ^{-}\circ \left( \alpha ^{+}\right)
			^{-1}(1)\right) }{\sup }\psi _{n}(\tau ).
	\end{eqnarray*}%
	So, 
	\begin{eqnarray}
	\left\Vert p(t)\right\Vert _{L^{2}(\alpha (t),\beta (t))}^{2} &\leq
	&\left\Vert \left( \widetilde{p},\widetilde{q}\right) \right\Vert
	_{L^{2}(0,1)}^{2}\underset{\tau \in \left[ 0,\alpha ^{-}\circ \left( \alpha
		^{+}\right) ^{-1}(1)\right) }{\sup }\psi _{n-1}(\tau )  \label{aa} \\
	&&+\left\Vert \left( \widetilde{p},\widetilde{q}\right) \right\Vert
	_{L^{2}(0,1)}^{2}\underset{\tau \in \left[ \alpha ^{-}\circ \left( \alpha
		^{+}\right) ^{-1}(1),\phi (0)\right) }{\sup }\psi _{n-1}(\tau )  \notag \\
	&&+\left\Vert \left( \widetilde{p},\widetilde{q}\right) \right\Vert
	_{L^{2}(0,1)}^{2}\underset{\tau \in \left[ 0,\alpha ^{-}\circ \left( \alpha
		^{+}\right) ^{-1}(1)\right) }{\sup }\psi _{n}(\tau ).  \notag
	\end{eqnarray}%
	\textbf{Case 2}: $t\in $ $\left[ \left( \alpha ^{-}\right) ^{-1}\circ \phi ^{%
		\left[ n\right] }(0),\left( \alpha ^{-}\right) ^{-1}\circ \phi ^{\left[ n%
		\right] }\circ \alpha ^{-}\circ \left( \alpha ^{+}\right) ^{-1}(1)\right) .$
	
	In this case, $p(t)$ might be expressed in function of $%
	p_{2n}(t),p_{2n+1}(t),p_{2n+2}(t)$ 
	\begin{equation}
	p(t,x)=\left\{ 
	\begin{array}{ccc}
	p_{2n}(t,x),\text{ \ } & \text{\textrm{if}} & x\in I_{4}(t):=\left[ t-\phi ^{%
		\left[ n\right] }(0),\beta (t)\right) ,\text{ \ \ \ \ \ \ \ \ \ \ \ \ \ \ \
		\ \ \ \ \ \ \ \ \ \ \ \ } \\ 
	p_{2n+1}(t,x), & \text{\textrm{if}} & x\in I_{5}(t):=\left[ t-\phi ^{\left[ n%
		\right] }\circ \alpha ^{-}\circ \left( \alpha ^{+}\right) ^{-1}(1),t-\phi ^{%
		\left[ n\right] }(0)\right) , \\ 
	p_{2n+2}(t,x), & \text{\textrm{if}} & x\in I_{6}(t):=\left[ \alpha
	(t),t-\phi ^{\left[ n\right] }\circ \alpha ^{-}\circ \left( \alpha
	^{+}\right) ^{-1}(1)\right) .\text{\ \ \ \ \ \ \ \ }%
	\end{array}%
	\right.  \label{p2}
	\end{equation}%
	In the same way, we obtain the estimate 
	\begin{eqnarray}
	\left\Vert p(t)\right\Vert _{L^{2}(\alpha (t),\beta (t))}^{2} &\leq
	&\left\Vert \left( \widetilde{p},\widetilde{q}\right) \right\Vert
	_{L^{2}(0,1)}^{2}\underset{\tau \in \left[ \alpha ^{-}\circ \left( \alpha
		^{+}\right) ^{-1}(1),\phi (0)\right) }{\sup }\psi _{n-1}(\tau )  \label{bb}
	\\
	&&+\left\Vert \left( \widetilde{p},\widetilde{q}\right) \right\Vert
	_{L^{2}(0,1)}^{2}\underset{\tau \in \left[ 0,\alpha ^{-}\circ \left( \alpha
		^{+}\right) ^{-1}(1)\right) }{\sup }\psi _{n}(\tau )  \notag \\
	&&+\left\Vert \left( \widetilde{p},\widetilde{q}\right) \right\Vert
	_{L^{2}(0,1)}^{2}\underset{\tau \in \left[ \alpha ^{-}\circ \left( \alpha
		^{+}\right) ^{-1}(1),\phi (0)\right) }{\sup }\psi _{n}(\tau ).  \notag
	\end{eqnarray}%
	Analogously, we have for the component $q$:\newline
	\textbf{Case 1}: $t\in \left[ \left( \beta ^{+}\right) ^{-1}\circ \xi ^{%
		\left[ n-1\right] }\circ \beta ^{+}\circ \left( \beta ^{-}\right)
	^{-1}(0),\left( \beta ^{+}\right) ^{-1}\circ \xi ^{\left[ n\right]
	}(1)\right) .$
	
	The expression of $q(t)$ might involve the expressions of $%
	q_{2n-1}(t),q_{2n}(t),q_{2n+1}(t)$ 
	\begin{equation}
	q(t,x)=\left\{ 
	\begin{array}{ccc}
	q_{2n-1}(t,x), & \text{\textrm{if}} & \text{ \ }x\in J_{1}(t):=\left[ \alpha
	(t),\xi ^{\left[ n-1\right] }\circ \beta ^{+}\circ \left( \beta ^{-}\right)
	^{-1}(0)-t\right) ,\text{ \ \ \ \ \ \ \ } \\ 
	q_{2n}(t,x),\text{ \ } & \text{\textrm{if}} & \text{ }x\in J_{2}(t):=\left[
	\xi ^{\left[ n-1\right] }\circ \beta ^{+}\circ \left( \beta ^{-}\right)
	^{-1}(0)-t,\xi ^{\left[ n\right] }(1)-t\right) , \\ 
	q_{2n+1}(t,x), & \text{\textrm{if}} & x\in J_{3}(t):=\left[ \xi ^{\left[ n%
		\right] }(1)-t,\beta (t)\right) .\text{ \ \ \ \ \ \ \ \ \ \ \ \ \ \ \ \ \ \
		\ \ \ \ \ \ \ \ \ \ \ }%
	\end{array}%
	\right.  \label{q1}
	\end{equation}%
	So, we have 
	\begin{eqnarray}
	&&\left\Vert q(t)\right\Vert _{L^{2}(\alpha (t),\beta (t))}^{2}  \label{ee}
	\\
	&\leq &\sum_{k=1}^{3}\int_{(t,x)\in \Sigma _{2n+k-2}^{q}}\left\vert
	q_{2n+k-2}(t,x)\right\vert ^{2}dx  \notag \\
	&\leq &\left\Vert \left( \widetilde{p},\widetilde{q}\right) \right\Vert
	_{L^{2}(0,1)}^{2}\times  \notag \\
	&&\sum_{k=1}^{3}\underset{x,(t,x)\in \Sigma _{2n+k-2}^{q}}{\sup }%
	\prod\limits_{i=0}^{n-2+\left[ \frac{k-1}{2}\right] }\left\vert F\left(
	\left( \alpha ^{-}\right) ^{-1}\circ \left( \phi ^{-1}\right) ^{[k]}\circ
	\beta ^{-}\circ \left( \beta ^{+}\right) ^{-1}(x+t)\right) \right\vert . 
	\notag
	\end{eqnarray}%
	By definition of the regions $\Sigma _{n}^{q},$ $n\geq 0$ given in (\ref%
	{region0})-(\ref{region1}), we have 
	\begin{equation}
	(t,x)\in \Sigma _{2n}^{q}\Leftrightarrow \text{\ }t+x\in \left[ \xi ^{\left[
		n-1\right] }\circ \beta ^{+}\circ \left( \beta ^{-}\right) ^{-1}(0),\xi ^{%
		\left[ n\right] }(1)\right) ,  \label{a}
	\end{equation}%
	\begin{equation}
	(t,x)\in \Sigma _{2n+1}^{q}\Leftrightarrow t+x\in \left[ \xi ^{\left[ n%
		\right] }(1),\xi ^{\left[ n\right] }\circ \beta ^{+}\circ \left( \beta
	^{-}\right) ^{-1}(0)\right) ,  \label{b}
	\end{equation}%
	and since $\xi $ is defined as 
	\begin{equation*}
	\xi =\beta ^{+}\circ \left( \beta ^{-}\right) ^{-1}\circ \phi \circ \beta
	^{-}\circ \left( \beta ^{+}\right) ^{-1},
	\end{equation*}%
	(\ref{a}) and (\ref{b}) turns to 
	\begin{eqnarray}
	(t,x) &\in &\Sigma _{2n}^{q}\Leftrightarrow  \label{c} \\
	t+x &\in &\left[ \beta ^{+}\circ \left( \beta ^{-}\right) ^{-1}\circ \phi ^{%
		\left[ n-1\right] }(0),\beta ^{+}\circ \left( \beta ^{-}\right) ^{-1}\circ
	\phi ^{\left[ n\right] }\circ \beta ^{-}\circ \left( \beta ^{+}\right)
	^{-1}(1)\right) ,  \notag
	\end{eqnarray}%
	\begin{eqnarray}
	(t,x) &\in &\Sigma _{2n+1}^{q}\Leftrightarrow  \label{d} \\
	t+x &\in &\left[ \beta ^{+}\circ \left( \beta ^{-}\right) ^{-1}\circ \phi ^{%
		\left[ n\right] }\circ \beta ^{-}\circ \left( \beta ^{+}\right)
	^{-1}(1),\beta ^{+}\circ \left( \beta ^{-}\right) ^{-1}\circ \phi ^{\left[ n%
		\right] }(0)\right) ,  \notag
	\end{eqnarray}%
	therefore, there exist $\chi _{1}:=\chi ^{n}(t,x)\in \left[ 0,\phi \circ
	\beta ^{-}\circ \left( \beta ^{+}\right) ^{-1}(1)\right) $ and $\chi
	_{2}:=\chi _{2}^{n}(t,x)\in \left[ \phi \circ \beta ^{-}\circ \left( \beta
	^{+}\right) ^{-1}(1),\phi (0)\right) $ such that 
	\begin{eqnarray*}
		(t,x) &\in &\Sigma _{2n}^{q}\text{ \ \ }\Leftrightarrow t+x=\beta ^{+}\circ
		\left( \beta ^{-}\right) ^{-1}\circ \phi ^{\left[ n-1\right] }\left( \chi
		_{1}\right) , \\
		(t,x) &\in &\Sigma _{2n+1}^{q}\Leftrightarrow t+x=\beta ^{+}\circ \left(
		\beta ^{-}\right) ^{-1}\circ \phi ^{\left[ n\right] }(\chi _{2}).
	\end{eqnarray*}%
	Thus, by combining (\ref{ee}),(\ref{c}) and (\ref{d}), we obtain
	
	\begin{eqnarray*}
		&&\sum_{k=1}^{3}\underset{x,(t,x)\in \Sigma _{2n+k-2}^{q}}{\sup }%
		\prod\limits_{i=0}^{n-2+\left[ \frac{k-1}{2}\right] }\left\vert F\left(
		\left( \alpha ^{-}\right) ^{-1}\circ \left( \phi ^{-1}\right) ^{[k]}\circ
		\beta ^{-}\circ \left( \beta ^{+}\right) ^{-1}(x+t)\right) \right\vert \\
		&\leq &\underset{\chi _{2}\in \left[ \phi \circ \beta ^{-}\circ \left( \beta
			^{+}\right) ^{-1}(1),\phi (0)\right) }{\sup }\prod\limits_{i=0}^{n-2}\left%
		\vert 
		\begin{array}{c}
			F\left( \alpha ^{-}\right) ^{-1}\circ \left( \phi ^{-1}\right) ^{[i]}\circ
			\beta ^{-}\circ \left( \beta ^{+}\right) ^{-1} \\ 
			\circ \beta ^{+}\circ \left( \beta ^{-}\right) ^{-1}\circ \phi ^{\left[ n-1%
				\right] }(\chi _{2})%
		\end{array}%
		\right\vert \\
		&&+\underset{\chi _{1}\in \left[ 0,\phi \circ \beta ^{-}\circ \left( \beta
			^{+}\right) ^{-1}(1)\right) }{\sup }\prod\limits_{i=0}^{n-2}\left\vert
		F\left( 
		\begin{array}{c}
			\left( \alpha ^{-}\right) ^{-1}\circ \left( \phi ^{-1}\right) ^{[i]}\circ
			\beta ^{-}\circ \left( \beta ^{+}\right) ^{-1} \\ 
			\circ \beta ^{+}\circ \left( \beta ^{-}\right) ^{-1}\circ \phi ^{\left[ n-1%
				\right] }\left( \chi _{1}\right)%
		\end{array}%
		\right) \right\vert \\
		&&+\underset{\chi _{2}\in \left[ \phi \circ \beta ^{-}\circ \left( \beta
			^{+}\right) ^{-1}(1),\phi (0)\right) }{\sup }\prod\limits_{i=0}^{n-1}\left%
		\vert F\left( 
		\begin{array}{c}
			\left( \alpha ^{-}\right) ^{-1}\circ \left( \phi ^{-1}\right) ^{[i]}\circ
			\beta ^{-}\circ \left( \beta ^{+}\right) ^{-1} \\ 
			\circ \beta ^{+}\circ \left( \beta ^{-}\right) ^{-1}\circ \phi ^{\left[ n%
				\right] }(\chi _{2})%
		\end{array}%
		\right) \right\vert \\
		&=&\underset{\chi _{2}\in \left[ \phi \circ \beta ^{-}\circ \left( \beta
			^{+}\right) ^{-1}(1),\phi (0)\right) }{\sup }\prod\limits_{i=0}^{n-2}\left%
		\vert F\left( \alpha ^{-}\right) ^{-1}\circ \phi ^{\left[ n-i-1\right]
		}(\chi _{2})\right\vert \\
		&=&\underset{\chi _{1}\in \left[ 0,\phi \circ \beta ^{-}\circ \left( \beta
			^{+}\right) ^{-1}(1)\right) }{\sup }\prod\limits_{i=0}^{n-2}\left\vert
		F\left( \left( \alpha ^{-}\right) ^{-1}\circ \phi ^{\left[ n-i-1\right]
		}\left( \chi _{1}\right) \right) \right\vert \\
		&=&\underset{\chi _{2}\in \left[ \phi \circ \beta ^{-}\circ \left( \beta
			^{+}\right) ^{-1}(1),\phi (0)\right) }{\sup }\prod\limits_{i=0}^{n-1}\left%
		\vert F\left( \left( \alpha ^{-}\right) ^{-1}\circ \phi ^{\lbrack n-i]}(\chi
		_{2})\right) \right\vert .
	\end{eqnarray*}%
	Finally, we get 
	\begin{eqnarray}
	\left\Vert q(t)\right\Vert _{L^{2}(\alpha (t),\beta (t))}^{2} &\leq
	&C\left\Vert \left( \widetilde{p},\widetilde{q}\right) \right\Vert
	_{L^{2}(0,1)}^{2}\underset{\chi \in \left[ \phi \circ \beta ^{-}\circ \left(
		\beta ^{+}\right) ^{-1}(1),\phi (0)\right) }{\sup }\psi _{n-1}(\chi )
	\label{cc} \\
	&&+C\left\Vert \left( \widetilde{p},\widetilde{q}\right) \right\Vert
	_{L^{2}(0,1)}^{2}\underset{\chi \in \left[ 0,\phi \circ \beta ^{-}\circ
		\left( \beta ^{+}\right) ^{-1}(1)\right) }{\sup }\psi _{n-1}(\chi )  \notag
	\\
	&&+C\left\Vert \left( \widetilde{p},\widetilde{q}\right) \right\Vert
	_{L^{2}(0,1)}^{2}\underset{\chi \in \left[ \phi \circ \beta ^{-}\circ \left(
		\beta ^{+}\right) ^{-1}(1),\phi (0)\right) }{\sup }\psi _{n}(\chi ).  \notag
	\end{eqnarray}%
	\textbf{Case 2}: $t\in \left[ \left( \beta ^{+}\right) ^{-1}\circ \xi ^{%
		\left[ n\right] }(1),\left( \beta ^{+}\right) ^{-1}\circ \xi ^{\left[ n%
		\right] }\circ \beta ^{+}\circ \left( \beta ^{-}\right) ^{-1}(0)\right) .$
	
	As previously, $q(t)$ might involve the values of $%
	q_{2n}(t),q_{2n+1}(t),q_{2n+2}(t)$ 
	\begin{equation}
	q(t,x)=\left\{ 
	\begin{array}{ccc}
	q_{2n}(t,x),\text{ \ } & \text{\textrm{if}} & x\in J_{4}(t):=\left[ \alpha
	(t),\xi ^{\left[ n\right] }(1)-t\right) ,\text{ \ \ \ \ \ \ \ \ \ \ \ \ \ \
		\ \ \ \ \ \ \ \ \ \ \ \ } \\ 
	q_{2n+1}(t,x), & \text{\textrm{if}} & x\in J_{5}(t):=\left[ \xi ^{\left[ n%
		\right] }(1)-t,\xi ^{\left[ n\right] }\circ \beta ^{+}\circ \left( \beta
	^{-}\right) ^{-1}(0)-t\right) , \\ 
	q_{2n+2}(t,x), & \text{\textrm{if}} & x\in J_{6}(t):=\left[ \xi ^{\left[ n%
		\right] }\circ \beta ^{+}\circ \left( \beta ^{-}\right) ^{-1}(0)-t,\beta
	(t)\right) .\text{\ \ \ \ \ \ \ \ }%
	\end{array}%
	\right.  \label{q2}
	\end{equation}%
	In the same way, the following estimate holds 
	\begin{eqnarray}
	\left\Vert q(t)\right\Vert _{L^{2}(\alpha (t),\beta (t))}^{2} &\leq
	&C\left\Vert \left( \widetilde{p},\widetilde{q}\right) \right\Vert
	_{L^{2}(0,1)}^{2}\underset{\chi \in \left[ 0,\phi \circ \beta ^{-}\circ
		\left( \beta ^{+}\right) ^{-1}(1)\right) }{\sup }\psi _{n-1}(\chi )
	\label{dd} \\
	&&+C\left\Vert \left( \widetilde{p},\widetilde{q}\right) \right\Vert
	_{L^{2}(0,1)}^{2}\underset{\chi \in \left[ \phi \circ \beta ^{-}\circ \left(
		\beta ^{+}\right) ^{-1}(1),\phi (0)\right) }{\sup }\psi _{n-1}(\chi )  \notag
	\\
	&&+C\left\Vert \left( \widetilde{p},\widetilde{q}\right) \right\Vert
	_{L^{2}(0,1)}^{2}\underset{\chi \in \left[ 0,\phi \circ \beta ^{-}\circ
		\left( \beta ^{+}\right) ^{-1}(1)\right) }{\sup }\psi _{n}(\chi ).  \notag
	\end{eqnarray}%
	From (\ref{aa}),(\ref{bb}),(\ref{cc}) and (\ref{dd}), we deduce that 
	\begin{equation}
	\underset{\tau \in \left[ 0,\phi (0)\right) }{\sup }\psi _{n}(\tau )\underset%
	{n\rightarrow \infty }{\longrightarrow }0\Longrightarrow \left\Vert
	(p,q)\right\Vert _{L^{2}(\alpha (t),\beta (t))}\underset{t\rightarrow \infty 
	}{\longrightarrow }0,  \label{n}
	\end{equation}%
	which finishes the proof of the first statement of Theorem \ref{theorem stab}%
	. The proof of the second and the third statements are just a consequences
	of (\ref{n}). By definition of the regions $\Sigma _{n}^{p},\Sigma _{n}^{q},$
	$n\geq 0,$ given in (\ref{region0})-(\ref{region1}), we can see that letting 
	$t\longrightarrow \infty $ is the same as $\phi ^{\left[ n\right] }(\tau
	)\longrightarrow \infty ,$ $\forall \tau \in \left[ 0,\phi (0)\right) ,$ so,
	if there exists a positive function $g$ such that 
	\begin{equation*}
	Cg\left( \phi ^{\left[ n\right] }(\tau )\right) \underset{n\rightarrow
		\infty }{\sim }\psi _{n}(\tau ),\text{ }\forall \tau \in \lbrack 0,\phi (0)),
	\end{equation*}%
	then obviously (\ref{rate}) holds. In particular, exponential stability
	follows immediately from 
	\begin{equation*}
	\underset{\tau \in \lbrack 0,\phi (0))}{\sup }\psi _{n}(\tau )=\underset{%
		\tau \in \lbrack 0,\phi (0))}{\sup }\exp \left[ \phi ^{\left[ n\right]
	}(\tau )\left( \frac{\ln \psi _{n}(\tau )}{\phi ^{\left[ n\right] }(\tau )}%
	\right) \right] .
	\end{equation*}
	
	The proof of the necessary part is straightforward. From (\ref{p2n+1}),(\ref%
	{p2n+2}) and (\ref{pppp}), we have
	
	\begin{eqnarray*}
		&&\int_{(t,x)\in \Sigma _{2n+1}^{p}}\left\vert p_{2n+1}(t,x)\right\vert
		^{2}dx+\int_{(t,x)\in \Sigma _{2n+1}^{p}}\left\vert p_{2n+1}(t,x)\right\vert
		^{2}dx \\
		&\geq &C\left\Vert \widetilde{q}\right\Vert _{L^{2}(0,1)}^{2}\underset{%
			x,(t,x)\in \Sigma _{2n+1}^{p}}{\inf }\prod\limits_{i=0}^{n}\left\vert
		F\left( \left( \alpha ^{-}\right) ^{-1}\circ \left( \phi ^{-1}\right)
		^{[i]}(t-x)\right) \right\vert \\
		&&+C\left\Vert p\right\Vert _{L^{2}(0,1)}^{2}\underset{x,(t,x)\in \Sigma
			_{2n+2}^{p}}{\inf }\prod\limits_{i=0}^{n}\left\vert F\left( \left( \alpha
		^{-}\right) ^{-1}\circ \left( \phi ^{-1}\right) ^{[i]}(t-x)\right)
		\right\vert \\
		&\geq &C\left( \left\Vert \widetilde{q}\right\Vert
		_{L^{2}(0,1)}^{2}+\left\Vert \widetilde{p}\right\Vert
		_{L^{2}(0,1)}^{2}\right) \times \\
		&&\left[ \underset{\tau \in \left[ 0,\alpha ^{-}\circ \left( \alpha
			^{+}\right) ^{-1}(1)\right) }{\inf }\psi _{n}(\tau )+\underset{\tau \in %
			\left[ \alpha ^{-}\circ \left( \alpha ^{+}\right) ^{-1}(1),\phi (0)\right) }{%
			\inf }\psi _{n}(\tau )\right] ,
	\end{eqnarray*}%
	therefore, if (\ref{assum}) is not satisfied then clearly stability cannot
	occur.
	
	Let us \ prove the second claim of Theorem \ref{theorem stab}. If $f\equiv 1$
	then $F\equiv 0.$ In this case, we infer from the exact formula of solutions
	given in\ (\ref{y1-}),(\ref{y1+}) and (\ref{y2+}) that we have $p_{1}\equiv
	0 $ while $q_{1}\neq 0$, and since $q$ is constant along the characteristic
	lines, $q$ is identically zero from the time that $q_{2}$ will be zero, that
	is $\ t\geq T^{\ast }=\left( \alpha ^{+}\right) ^{-1}\circ \beta ^{+}\circ
	\left( \beta ^{-}\right) ^{-1}(0)$ which is the same time for boundary
	controllability of system (\ref{control}).
	
	\section{Further remarks and open questions}
	
	Let us discuss briefly some possible variations and generalization of the
	obtained results in this work.
	
	\begin{itemize}
		\item It is our hope that the tools developed in this paper may help in
		dealing with the distributed control case 
		\begin{equation}
		\left\{ 
		\begin{array}{ccc}
		y_{tt}(t,x)=y_{xx}(t,x)+\chi _{\omega _{T}}h(t,x), & \mathrm{in} & Q_{T},%
		\text{ \ } \\ 
		y(t,\alpha (t))=\text{ }y(t,\beta (t))=0,\text{ \ \ \ \ \ \ \ } & \mathrm{in}
		& (0,T), \\ 
		y(0,x)=y_{0}(x),\text{ }y_{t}(0,x)=y_{1}(x), & \mathrm{in} & (0,1),%
		\end{array}%
		\right.  \label{A}
		\end{equation}
		
		where $\omega _{T}$ is a moving subset of $Q_{T}:=(0,T)\times (0,1)$ defined
		by 
		\begin{equation*}
		\omega _{T}=\left\{ (t,x)\in Q_{T},\text{ }x\in (a(t),b(t))\right\} ,
		\end{equation*}%
		and $\left( y_{0},y_{1},h\right) \in H_{0}^{1}(0,1)\times L^{2}(0,1)\times
		L^{2}(\omega _{T})$. Actually, we can determine the minimal time for which
		the time-dependent geometric control condition introduced in \cite[
		Definition 1.6]{Le Rousseau} is satisfied. The latter condition states that
		every generalized bicharacteristic must meet the moving control region at
		some time $T.$ It is easy to verify this condition in the one dimensional
		settings. Indeed, under assumption (\ref{assumption}) with $a,b\in
		C^{1}(0,T) $ and $\left\Vert a^{\prime }\right\Vert _{L^{\infty
			}(0,T)},\left\Vert b^{\prime }\right\Vert _{L^{\infty }(0,T)}<1,$ we find
		that all the characteristics with positive slope or negative slope emerging
		from the point $(0,x),$ for any $x\in (0,1)$ meet $\omega _{T}$ if, and only
		if $T>T^{\ast }$ where $T^{\ast }$ is given by 
		\begin{equation*}
		T^{\ast }=\max \left\{ T_{1},T_{2}\right\} =\max \left\{ b^{+}\circ \beta
		^{+}\circ \left( \beta ^{-}\right) ^{-1}\circ b(0),a^{-}\circ \alpha
		^{-}\circ \left( \alpha ^{+}\right) ^{-1}\circ a(0)\right\} .
		\end{equation*}%
		In particular, if $\alpha \equiv 0$ and $\beta \equiv 1,$ the time $T^{\ast
		} $ is given by $T^{\ast }=2\max \left\{ a,1-b\right\} $ which is exactly
		the time given in \cite{Zuazua}.
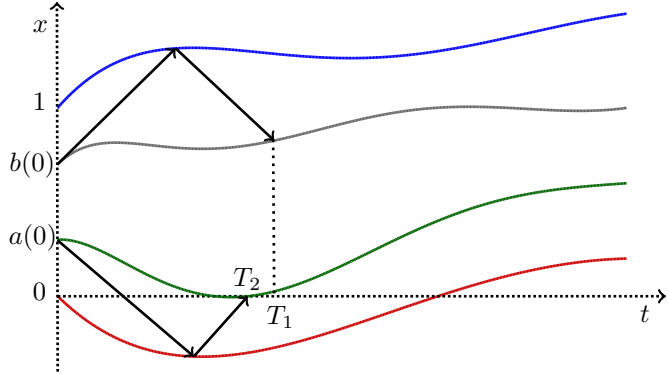
		\begin{figure}[H]
			\centering
			\definecolor{qqwuqq}{rgb}{0,0.39215686274509803,0}  %
			\definecolor{wwwwww}{rgb}{0.4,0.4,0.4} \definecolor{ccqqqq}{rgb}{0.8,0,0}  %
			\definecolor{qqqqff}{rgb}{0,0,1}  

		\end{figure}

		Controllability and stabilizability of the multidimensional wave equation in
		non-cylindrical domains has been investigated by Bardos and Chen in \cite%
		{Bardos}. By assuming that the domain is expending, exact controllability and
		stability have been established by a control and a frictional
		damping acting on the entire domain.
		
		Distributed controllability of system (\ref{A}) has been studied in \cite%
		{Castro} in a cylindircal domain with moving control support, i.e. $\alpha
		\equiv 0$ and $\beta \equiv 1$. It has been proved that exact
		controllability holds if the moving control support $\omega _{T}$ satisfies
		the geometric control condition without the restriction $\left\Vert
		a^{\prime }\right\Vert _{L^{\infty }(0,T)},\left\Vert b^{\prime }\right\Vert
		_{L^{\infty }(0,T)}<1$. It worths to mention that problem (\ref{A}) has been
		recently studied in \cite{Cui3} with a very partiuclar boundary curves and
		moving control support. The critical time of control seems to be far from
		being optimal.
		
		\item Note that we have not used the $L^{2}$ settings in a crucial way. The
		same results can be proved in $L^{p},$ $p\in \lbrack 1,\infty ),$ or in the
		space of continuous functions.
	\end{itemize}
	
	\section*{Acknowledgment}
	
	I would like to thank my supervisor, Ammar Khodja Farid, for his continuous
	support and valuable remarks, as well the anonymous referees for their
	comments and suggestions.

\end{document}